\documentclass{article}
\usepackage{amsmath, amssymb, amsthm, graphicx,bm,
cite,
comment}
\usepackage[all]{xy}
\usepackage{titlesec}
\DeclareMathOperator{\gr}{Graph}
\DeclareMathOperator{\spt}{Supp}
\DeclareMathOperator{\re}{Re}
\theoremstyle{plain}% default
\newtheorem{thm}{Theorem}[section]
\newtheorem{lem}[thm]{Lemma}
\newtheorem{cor}[thm]{Corollary}

\newtheorem*{thm*}{Theorem}
\newtheorem*{lem*}{Lemma}
\newtheorem*{cor*}{Corollary}

\theoremstyle{definition}
\newtheorem{dfn}[thm]{Definition}

\newtheorem{hp}[thm]{Hypothesis}

\newtheorem*{dfn*}{Definition}

\newtheorem*{hp*}{Hypothesis}

\numberwithin{equation}{section}\numberwithin{figure}{section}

\def\e#1\e{\begin{equation}#1\end{equation}}
\def\iz#1\iz{\begin{itemize}#1\end{itemize}}
\def\ea#1\ea{\e{\begin{split}#1\end{split}}\e}
\def\eq{\eqref}
\def\l{\label}
\def\0{\hspace{0pt}}

\def\ph{\phi}

\def\Uh_#1{\,\widehat{\!U}_{\!#1}}

\def\Ps{\Psi}
\def\Ph{\Phi}
\def\r{\rho}

\def\et{\eta}
\def\ps{\psi}
\def\ch{\chi}
\def\ft{\mathfrak t}

\def\fsu{\mathfrak{su}}
\def\px{\approx}

\def\Im{\mathop{\rm Im}}

\def\ge{\geqslant}
\def\le{\leqslant\nobreak}
\def\les{\lesssim\nobreak}

\def\cB{{\mathbin{\cal B}}}

\def\cF{{\mathbin{\cal F}}}
\def\cG{{\mathbin{\cal G}}}
\def\cH{{\mathbin{\cal H}}}

\def\cJ{{\mathbin{\cal J}}}
\def\cK{{\mathbin{\cal K}}}
\def\cL{{\mathbin{\cal L}}}
\def\cM{{\mathbin{\cal M}}}

\def\O{{\mathbin{\cal O}}}
\def\cP{{\mathbin{\cal P}}}
\def\cR{{\mathbin{\cal R}}}

\def\cU{{\mathbin{\cal U}}}

\def\={\equiv}

\def\cT{{\mathbin{\cal{T}}}}
\def\cQ{{\mathbin{\cal Q}}}
\def\cW{{\mathbin{\cal W}}}

\def\C{{{\mathbb C}}}

\def\CP{{{\C}P}}
\def\RP{{{\R}P}}

\def\R{{\mathbin{\mathbb R}}}
\def\Z{{\mathbin{\mathbb Z}}}

\def\al{\alpha}
\def\be{\beta}
\def\ga{\gamma}
\def\de{\delta}
\def\io{\iota}
\def\ep{\epsilon}
\def\la{\lambda}

\def\th{\theta}
\def\ta{{\tau}}
\def\ze{\zeta}

\def\si{\sigma}
\def\om{\omega}
\def\De{\Delta}

\def\Si{\Sigma}

\def\Th{\Theta}
\def\Om{\Omega}
\def\Ga{\Gamma}
\def\Up{\Upsilon}

\def\d{\partial}
\def\db{\,\ov{\!\partial}}
\def\ts{\textstyle}

\def\w{\wedge}
\def\-{\setminus}
\def\bu{\bullet}
\def\op{\oplus}

\def\ot{\otimes}
\def\bot{\bigotimes}
\def\ov{\overline}

\def\ul{\underline}
\def\iy{\infty}

\def\t{\times}

\def\nb{\nabla}
\def\sb{\subset}

\begin{document}
\date{}
\title{Example of Compact Special Lagrangians with a Stable Singularity}
\author{Yohsuke Imagi}
\maketitle
%This paper contains the first example in the literature of compact special Lagrangians with isolated conical singularities to which one can make the full application of an established local moduli-theory.
%\tableofcontents
\begin{abstract}
We construct a family of compact almost Calabi--Yau manifolds of complex dimension 3 and therein a corresponding family of compact special Lagrangians with one-point singularities modelled upon that $T^2$-cone constructed by Harvey and Lawson \cite[Chapter III.3.A, Theorem 3.1]{HL} and characterized by Haskins \cite[Theorem A]{H} as a stable $T^2$-cone in the terminology by Joyce \cite[Definition 3.4 and Example 3.5]{J5}.
%We describe, moreover, how the general results by Joyce \cite[\S10]{J5} and the present author \cite{I} apply to our circumstances.
\end{abstract}

\section*{Introduction}
This paper contains two sections: \S\ref{con} is devoted to an abstract gluing construction of compact special Lagrangians with one-point singularities modelled upon Harvey--Lawson's $T^2$-cone; and \S\ref{ex} applies it to a concrete example.

The building blocks of the gluing construction contain a nodal complex 3-fold $N_0$ with an almost Calabi--Yau structure in the sense of Definition \ref{N}; a compact special Lagrangian $X_0$ in $N_0$ with one-point singularity $n_0,$ which is a node of $N_0,$ modelled upon another $T^2$-cone $C;$ and a smooth prime divisor $D_0$ on $N_0$ which intersects $X_0$ only at $n_0.$
In the example of \S\ref{ex}, $N_0$ is a quintic hypersurface in $\CP^4,$ equipped with the Fubini--Study metric; $X_0$ is the fixed-point set of an antiholomorphic involution of $N_0;$ and $D_0$ is a projective plane in $N_0.$

We denote by $N_1\to N_0$ the blow-up of $D_0.$
%We take on $N_1$ a K\"ahler metric in the sense of \S\ref{N} with respect to which the area of the exceptional fibre $\pi^{-1}(n_0)$ is small to a certain extent.
The normal bundle in $N_1$ to the exceptional fibre over $n_0$ admits a noncompact special Lagrangian $Y$ which is, near infinity, contained in $C$ and has one-point singularity modelled upon Harvey--Lawson's $T^2$-cone.
Gluing together $X_0$ and $Y$ we construct in $N_1$ a compact special Lagrangian $X_1$ with one-point singularity modelled upon Harvey--Lawson's $T^2$-cone.

%Since $Y$ is, near infinity, contained in $C$ it follows that the gluing process is not topologically so obstructed as in the fourth \cite{J4} of the five papers by Joyce \cite{J1,J2,J3,J4,J5}.
%The analytical obstruction also vanishes as in the second \cite{J2} of these five, because of the stability of Harvey--Lawson's $T^2$-cone.
%The remaining analysis is close to that by Pacini \cite{P3} though the latter circumstances are simpler in the respect that the ambient spaces are flat.

In the example of \S\ref{ex}, the complex $3$-fold $N_1,$ though strictly nodal, yet admits both smoothings and projective small resolutions; the special Lagrangian $X_1,$ homeomorphic to $S^1\t D^2$ with boundary collapsed to a point, persists under the smoothings of $N_1$ whereas necessarily desingularized into $S^3$ under the projective small resolutions of $N_1;$ and there are also in the smoothings of $N_1$ some desingularizations of $X_1$ which are, however, diffeomorphic to $S^1\t S^2.$

{\bf Acknowledgements.}
The general idea of applying a simultaneous gluing method to (almost) Calabi--Yau manifolds and to special Lagrangian submanifolds of them goes back to the works by Chan \cite{Ch1,Ch2,Ch3,Ch4}. Dominic Joyce suggested to me in March 2012 the possibility of applying this idea to the problem of finding compact special Lagrangians with stable $T^2$-cone singularities; on which thereafter I had more detailed discussions with Mark Haskins, especially whilst visiting Imperial College London in November to December of 2014; during which, moreover, Johannes Nordstr\"om gave me the idea of finding the example of the quintic hypersurface $N_0\sb\CP^4.$
I had other helpful discussions about \S\ref{con} with Ryushi Goto, Kota Hattori, Shinichiroh Matsuo; and about \S\ref{ex} with Kei Irie and Toru Yoshiyasu on the topology of the special Lagrangians.
There were also helpful comments by a referee about the presentation of the paper.

Hein and Sun \cite{HS} proved in the meanwhile that every nodal quintic hypersurface in $\CP^4$ admits a truly Calabi--Yau (i.e.\ Ricci-flat) metric with a suitable asymptotic behaviour; which one can use perhaps for a truly Calabi--Yau version of our result.

The present paper was prepared during my work at Kyoto University, the Kavli Institute for Physics and Mathematics of the Universe, the Chinese University of Hong Kong, and Nagoya University; with financial support by the World Premier International Research Center Initiative, and the grants-in-aid (10J00699, 16K17587, 18J00075) of the Japan Society for the Promotion of Science.

\section{Gluing Construction}\l{con}
This section contains three subsections: \S\ref{11} describes the building blocks of the gluing construction; \S\ref{12} constructs an approximate solution of the gluing problem; and in \S\ref{13} we perturb it to the true solution.

\subsection{Building Blocks}\l{11}
We denote by $(x_1,x_2,x_3,x_4)$ the co\"ordinate-system of $\C^4,$ and by $Q_0$ the quadric in $\C^4$ defined by $x_1x_4-x_2x_3=0.$
We introduce a notion of almost Calabi--Yau structures on nodal complex $3$-folds, extending that by Joyce \cite[Definitions 2.8 and 2.9]{J5}:
\begin{dfn}\l{N}
By a {\it complex $3$-fold} we mean a complex $3$-dimensional analytic space.
We say that a complex 3-fold $N$ is {\it nodal} if for any singular point $n\in N$ there exists a biholomorphism-of-germs $(Q_0,0)\cong(N,n).$

Let $N$ be a nodal complex 3-fold; then, a {\it node} of $N$ shall mean a singular point of $N;$ and $N'$ shall denote the complement of $N$ to its nodes.
Moreover, a {\it K\"ahler form on $N$} shall mean a K\"ahler form $\om$ on $N'$ such that for any node $n\in N$ there exists a biholomorphism-of-germs $\Ph_n:(Q_0,0)\cong(N,n)$ with $\Ph_n^*\om$ extendable to a K\"ahler form near 0 in $\C^4;$ a {\it holomorphic volume-form} on $N$ an $\O_N$-module isomorphism $\Om_N^3\cong\O_N;$
and an {\it almost Calabi--Yau structure} on $N$ the pair $(\om,\Om)$ where $\om$ is a K\"ahler form, and $\Om$ a holomorphic volume-form, both on $N.$

Let $(\om,\Om)$ be an almost Calabi--Yau structure on $N;$ then, a {\it smooth special Lagrangian} in $(N;\om,\Om)$ shall mean a submanifold of $N'$ of real dimension 3 on which $\om$ and $\Im\Om$ vanish.
\end{dfn}

We put $\th_0:=i/2\sum_{j=1}^4dx_j\w d\bar{x}_j.$
We denote by $\Th_0$ the residue on $Q_0$ of the meromorphic 4-form $(x_1x_4-x_2x_3)^{-1}dx_1\w dx_2\w dx_3\w dx_4.$
% or more explicitly, $\Th_0=dx_2\w dx_3\w dx_4/x_4$ wherever $x_4\ne0$ for instance.
This $(\th_0,\Th_0)$ defines on $Q_0$ an almost Calabi--Yau structure.
It is invariant under the $T^2$-action on $Q_0$ defined by
\e\l{T^2}(x_1,x_2,x_3,x_4)\mapsto(t_1x_1,t_2^{-1}x_2,t_2x_3,t_1^{-1}x_4)\e
where $(t_1,t_2)\in T^2\sb(\C^*)^2.$

%In general, the fixed-point set of an antisymplectic antiholomorphic involution of an almost Calabi--Yau manifold is unless empty a smooth special Lagrangian; in particular, 
%We put $C:=\{(x_1,x_2,\bar{x}_2,\bar{x}_1)\in Q_0:|x_1|=|x_2|\},$ 
We denote by $C$ the fixed-point set in $Q_0$ of the involution $(x_1,x_2,x_3,x_4)\mapsto(\bar{x}_4,\bar{x}_3,\bar{x}_2,\bar{x}_1),$ which is topologically a $T^2$-cone; and $C':=C\-\{0\}$ is a $T^2$-invariant smooth special Lagrangian in $(Q_0;\th_0,\Th_0).$

We denote by $P_0$ the $T^2$-invariant complex plane in $Q_0$ defined by $x_3=x_4=0.$
We denote by $\pi:Q_1\to Q_0$ the blow-up of $P_0,$ which is a small resolution. We identify $Q_1$ with
\[\O_{\CP^1}(-1)^{\op2}=\{(x_1,x_2,x_3,x_4;y_1,y_2)\in\C^4\t\CP^1:x_1y_2=x_2y_1, x_3y_2=x_4y_1\};\]
and identify $\CP^1$ with the exceptional fibre of $\pi:Q_1\to Q_0.$

The $T^2$-action \eq{T^2} lifts to $Q_1.$
There is on $Q_1$ a $T^2$-invariant holomorphic volume-form $\Th_1:=\pi^*\Th_0.$
%$\pi^*\Om_Q_0=-dx_2\w dx_4 \w d(y_2/y_1)$
%wherever $y_1\ne0$ for instance.

We construct on $Q_1$ a $T^2$-invariant K\"ahler form.
We denote by $P_1$ the strict transform of $P_0$ under $\pi.$
There is over $Q_1$ the line-bundle $[-P_1],$ isomorphic to the pull-back of $\O_{\CP^1}(1)$ under the projection $Q_1\to\CP^1;$ and accordingly, on $[-P_1],$ a $T^2$-invariant Hermitian metric $h_{\rm FS}$ whose curvature is the pull-back of the Fubini--Study form on $\CP^1.$
% on  with $h_{{P_1}}=1+|y_2/y_1|^2$ wherever $y_1\ne0$ and with $h_{{P_1}}=|y_1/y_2|^2+1$ wherever $y_2\ne0.$
Let $h_{\rm flat}$ be a flat metric on $[-P_1]$ restricted to $Q_1\-{P_1};$ which exists because $[-P_1]$ is, by definition, trivial over $Q_1\-P_1.$
Let $\ch:[0,\iy)\to[0,1]$ be a cut-off function with $\ch\=1$ on $[0,1]$ and $\ch\=0$ on $[1,2].$
We make the identification $\ch=\ch(|x_3|^2+|x_4|^2);$ the latter is a $T^2$-invariant function on $Q_1.$
%Let $d_{P_1}:Q_1\to\R$ be the distance from ${P_1}$ relative to the K\"ahler form $\th_0+\om_{\CP^1}$ under the embedding $Q_1\sb\C^4\t\CP^1.$
%Let $\ch:\R_{\ge0}\to\R$ be a compactly-supported smooth function with $\ch\=1$ on $[0,1]$ and with $\ch\=0$ on $[2,\iy).$
Then, with $\de>0$ and small enough, there is on $Q_1$ a $T^2$-invariant K\"ahler form
\e\l{omti}
\th_1:=\pi^*\th_0+\de(i/2)\d\db\log[\ch h_{\rm FS}+(1-\ch)h_{\rm flat}].\e

The following theorem resembles the result by Ionel and Min-Oo \cite[Theorem 6.1]{IM} who use, however, a truly Calabi--Yau metric on $Q_1$:
\begin{thm}\l{NCBB}
With $\de$ made smaller if need be, there exists in $(Q_1;\th_1,\Th_1)$ a closed $T^2$-invariant special Lagrangian $Y$ such that $Y$ is, near infinity, contained in $\pi^*C';$ $Y$ has one-point singularity, at $p_Y:=(0,0,0,0;0,1)\in Q_1,$ modelled upon Harvey--Lawson's $T^2$-cone; and $Y':=Y\-\{p_Y\}$ is diffeomorphic to $\R\t T^2.$
\end{thm}
\begin{proof}
We denote by $U_s,$ $s>0,$ the open subset of $Q_1$ defined by $|x_3|^2+|x_4|^2<s.$
There is a unique $T^2$-moment map
\[m=(m_1,m_2):Q_1\to\ft^2\=\R^2\]
which satisfies, on $U_1,$
\[2m_1=|x_1|^2-|x_4|^2+\de\frac{|y_1|^2}{|y_1|^2+|y_2|^2}\]
and, on $Q_1\-U_2,$ $2m_1=|x_1|^2-|x_4|^2;$ and the same formulae for $m_2$ but with $-|x_2|^2+|x_3|^2$ in place of $|x_1|^2-|x_4|^2.$

Define $Y\sb Q_1$ by $m=0$ and $x_1x_4=x_2x_3\ge0.$
Then, $Y\-U_2=(\pi^*C)\-U_2;$ and moreover, on $Y\cap U_1,$
\e\l{YU}y_2\ne0\text{ and }|x_2|^2=|x_4|^2=\frac{\de|y|^2}{1-|y|^4}
\text{ where }y:=y_1/y_2.\e
This implies in turn that $Y\cap U_1$ is singular only at $p_Y,$ modelled upon Harvey--Lawson's $T^2$-cone
\[\{(x_2,x_4,y)\in\C^3:|x_2|^2=|x_4|^2=\de|y|^2,x_2x_4y\in\R_{\ge0}\}.\]
% relative to the metric $|dz|^2+|dx_4|^2+\de|dy|^2.$
The formula \eq{YU} implies moreover that $Y\cap U_1\cap\CP^1=\{p_Y\},$ that $\re x_2x_3>0$ on $Y'\cap U_1,$ and that $\re x_2x_3$ gives $Y'\cap U_1$ the structure of a principal $T^2$-bundle.

We prove that if $\de$ is small enough then the interior of $Y\-U_{1/2}$ is nonsingular.
We denote by $V$ the open subset of $Q_1$ defined by $\re x_2x_3>0.$
The definition of $C$ implies $\pi^*C\-U_{1/2}\sb V.$
Hence, recalling that if $\de$ were 0 then $Y$ would agree with $\pi^*C,$ one finds that if $\de$ is small enough then $Y\-U_{1/2}$ is so close to $\pi^*C\-U_{1/2}$ that $Y\-U_{1/2}\sb V.$
Accordingly, $Y\-U_{1/2}$ is outside $\CP^1;$ and lies in the zero set of the map $(m,\Im z_2,z_3):V\to\R^3$ which is a $\de$-perturbation of the other map $V\to\R^3$ defined by
\[(x_1,x_2,x_3,x_4)\mapsto(|x_1|^2-|x_4|^2,-|x_2|^2+|x_3|^2,\Im x_2x_3).\]
Computation implies that 0 is a regular value of the latter map and accordingly of its $\de$-perturbation with $\de$ small enough.
Consequently, the interior of $Y\-U_{1/2}$ is nonsingular.
% not intersecting $\pi^{-1}(0);$ and

Thus, $Y'$ is nonsingular; and $Y\cap\CP^1=\{p_Y\}.$
The arguments in the previous paragraph imply moreover that $\re x_2x_3>0$ on $Y'$ and that $\re x_2x_3:Y'\to\R_+$ is a principal $T^2$-bundle. 

Since the $T^2$-orbits in $Y'$ are of real codimension one it follows that $Y'$ is Lagrangian in $(Q_1,\th_1).$
We prove that $Y'$ is special Lagrangian with respect to $\Th_1.$
The $T^2$-action \eq{T^2} induces on $Q_0'$ the two vector-fields
\[v_1:=i(x_1\d_1-\bar{x}_1\ov{\d_1}-x_4\d_4+\bar{x}_4\ov{\d_4})\text{ and }
v_2:=i(-x_2\d_2+\bar{x}_2\ov{\d_2}+x_3\d_3-\bar{x}_3\ov{\d_3}).\]
If $x_4\ne0,$ for instance, then $\Th_1=dx_2\w dx_3\w dx_4/x_4;$ and accordingly,
\[v_2\lrcorner v_1\lrcorner\,\Im\Th_1=-\Im(x_2\,dx_3+x_3\,dx_2)=-d\,(\Im x_2x_3).\]
This identity extends by continuity to the whole $Q_1;$ and then, implies $\Im\Th_1|_{Y'}=0.$
\end{proof}

We extend the notion by Joyce \cite[Definition 3.6]{J5} of isolated conical singularities.
Define $r_C:C'\to\R_+$ by $x\mapsto |x|,$ i.e.\ the Euclidean distance from 0.
The {\it cylindrical} metric on $C'$ shall mean $|dx|^2/|x|^2,$ i.e.\ the conformal transform by $r_C^{-2}$ of the Euclidean metric restricted to $C'.$
Let $X$ be a subset of $Q_0$ which passes though the origin $0\in Q_0\sb\C^4;$
then, we say that $X$ {\it approaches $C$ with order $\la>2$} if there exists a homeomorphism-of-germs $\ph:(C,0)\to(X,0)$ which is, outside 0, diffeomorphic with $|\ph(x)-x|_{C^k}=O(|x|^{\la-1}),$ $k\in\Z_{\ge0},$ relative to the cylindrical metric on $C'.$

We make the following hypothesis, which holds in the example of \S\ref{ex}:
\begin{hp}\l{H}
There be a nodal complex 3-fold $N_0$ with an almost Calabi--Yau structure $(\om_0,\Om_0);$ therein a compact special Lagrangian $X_0$ with one-point singularity at a node $n_0\in N_0;$ on $N_0$ a smooth prime divisor $D_0$ which intersects $X_0$ only at $n_0;$ and a biholomorphism-of-germs $\Ph:(Q_0,0)\to(N_0,n_0)$ such that $\Ph^*\om_0|_0=\th_0,$ $\Ph^*\Om_0\=\Th_0$ near $0$ in $Q_0,$ $\Ph^*X_0$ approaches $C$ with order $>2,$ and $\Ph^*D_0\sb P_0.$
\end{hp}

We denote by $\pi: N_1\to N_0$ the blow-up of $D_0,$ abusing the same symbol as $\pi: Q_1\to Q_0.$
This $\pi: N_1\to N_0$ is a small morphism; in particular, $\Om_1:=\pi^*\Om_0$ is a holomorphic volume-form on $N_1.$

We construct on $N_1$ a K\"ahler form in the sense of Definition \ref{N}.
Let $D_1$ denote the strict transform of $D_0$ under $\pi: N_1\to N_0.$
Identify $\pi: N_1\to N_0$ locally with $\pi: Q_1\to Q_0$ by using $\Ph:(Q_0,0)\cong(N_0,n_0)$ and, for every other node $n\in D_0,$ any biholomorphism-of-germs $(Q_0,0)\cong(N_0,n).$
Let $h_1$ be a Hermitian metric near $D_1$ on the line-bundle $[-{D_1}]$ which agrees, near every exceptional fibre, with $h_{\rm FS}.$
Let $h_{\rm flat}$ be a flat metric on $[-D_1]$ over $N_1\-{D_1}.$
%Let $g_0$ be a Hermitian metric on $ N_1$ which, near every exceptional fibre in $ N_1,$ corresponds to the K\"ahler form $\th_0+\om_{\CP^1};$ and let $d_{D_1}: N_1\to\R_{\ge0}$ be the distance from ${D_1}$ relative to $g_0.$
Take a Hermitian metric on $N_1$ with respect to which the squared distance from $D_1$ agrees, near the exceptional fibre over $n_0,$ with $|x_3|^2+|x_4|^2:Q_1\to\R_{\ge0}.$ This distance we denote by $d:N_1\to\R_{\ge0}.$
Let $\ch:\R_{\ge0}\to\R_{\ge0}$ be that cut-off function used in \eq{omti}; and put $\ch_\ep:=\ch(d^2/\ep^2),$ $\ep>0.$
% and $\wh\ch:N_1\to[0,1]$ a cut-off function with 
%\[\wh\ch\=\begin{cases}1\text{ near }D_1, \\0\text{ outside a larger neighbourhood of $D_1,$ and}\\
%\ch\text{ near every exceptional fibre}.\end{cases}\]
%Let $\ep>0;$ which we identify with the fibrewise dilation of the normal bundle to $D_1$ in $N_1,$ and of $Q_1=\O_{\CP^1}(-1)^{\op2}.$
Then:
\begin{thm*}
With $\de>0$ and small enough, there is, associated to every $\ep\in(0,1],$ a K\"ahler form on $N_1$ of the form
\e\l{omt}\om_1:=\pi^*\om_0+\de\ep^2(i/2)\d\db\log[\ch_\ep h_1+(1-\ch_\ep)h_{\rm flat}].\e
\end{thm*}
\begin{proof}
It suffices to prove that $\om_1$ is positive definite on the annular region where $\ep_*\ch$ varies nontrivially.
Taking a connexion on the normal bundle to $D_1$ in $N_1$ one gets on this annular region a decomposition $\pi^*\om_0=\ph_h+\ph_v$ into the horizontal direction and the vertical direction.
One has then, by \eq{omt},
\[\ep^*\om_1=\ph_h+\ep^2(\ph_v+\de\ps_\ep)\]
where $\ps_\ep$ is $\ep$-uniformly bounded.
Consequently, with $\de$ small enough, $\om_1$ is positive definite on the annular region.
\end{proof}

\subsection{Construction of an Approximate Solution}\l{12}
Gluing $X_0,Y$ we construct in $(N_1;\om_1,\Om_1)$ a compact Lagrangian $Z$ with one-point singularity modelled upon Harvey--Lawson $T^2$-cone.
Since $Y$ is, near infinity, contained in $\pi^*C'$ it follows that this process is not so obstructed as in that by Joyce \cite{J4}.
%This process is not so obstructed as in \S7.3 and \S7.4 of Joyce \cite{J5}; for $Y$ is, near infinity, contained in $\pi^*C'.$

Let $K$ be a sufficiently-large precompact neighbourhood in $Q_1$ of its exceptional fibre, and identify it with the corresponding neighbourhood in $N_1$ of the excpetional fibre over $n_0.$
We denote by $\ep,$ for brevity's sake, the fibrewise $\ep$-dilation of the vector bundle $Q_1=\O_{\CP^1}(-1)^{\op2}.$
Then, by \eq{omt}, $\ep^{-2}\ep^*\om_1$ defines on $K$ a K\"ahler form which is, moreover, smooth at $\ep=0$ at which it is equal to $\th_1.$
One has also, pointwise on $\pi^{-1}(0),$ $\ep^{-2}\ep^*\om_1=\th_1.$
There exists, by Moser's technique, a Hamiltonian flow $\Ps_\ep:K\to Q_1$ with
$\ep^{-2}\Ps_\ep^*\ep^*\om_1=\th_1$
where $\Ps_\ep$ depends smoothly upon $\ep$ near and at 0; $\Ps_0$ is the identity-map of $Q_1;$ and $\Ps_\ep$ is, pointwise on $\pi^{-1}(0),$ the identity up to the linear order.
%Making $K$ larger if need be, one may suppose that $Y\-K=(\pi^*C)\-K$ and that $Y,\pi^*C$ agree also near the boundary of $K.$
%Making $K$ yet larger if need be, one may suppose $(\pi^*C)\-K$ so far from ${P_1}$ that $\om_1=\pi^*\om_0$ near $(\pi^*C)\-K.$

We put, for every $s>0,$ $B_s:=\{x\in Q_0\sb\C^4:|x|<s\}.$
% where $|\bu|$ denotes the Euclidean norm on $\C^4.$
The method by Chan \cite[Proof of Theorem 4.9]{Ch1} yields a constant $\r>0$ and an embedding $\Up:(B_\r,0)\to (Q_0,0)$ which is, outside 0, smooth with $\Up^*\om_0=\th_0$ and $|\Up(x)-x|_{C^k}=O(|x|^2),$ $k\in\Z_{\ge0},$ relative to the cylindrical metric $|dx|^2/|x|^2.$
The method by Joyce \cite[\S\S4.1 and 4.2]{J1} yields a dilation-invariant Weinstein neighbourhood of $C'$ in $(Q_0',\th_0),$ a constant $\la\in(2,3),$ and a smooth function $\xi_0:C'\cap B_\r\to\R$ such that
$\Up^*X_0=\gr d\xi_0$ with $|\xi_0|_{C^k}=O(r_C^\la),$ $k\in\Z_{\ge0},$ with respect to the cylindrical metric on $C'.$

%Since $\Ps_0$ is an identity and since $Y,\pi^*C$ agree near the boundary of $K$ one gets an $\ep$-independent constant $R>0$ with $\Up^*\pi_*\dl_{\ep*}\Ps_{\ep*}(Y\cap K)\-B_{\ep R}$ contained in the graph of $df_\ep$ where $f_\ep$ is a smooth $\R$-valued function defined in $C'$ outside $B_{\ep R}$ but near the boundary of $B_{\ep R}.$

There is now a compact Lagrangian in $(N_1,\om_1)$ of the form
\e\l{Z}Z:=\ep\Ps_\ep(Y\cap K)\cup(C\cap B_\r\-B_{\ep R})\sqcup(X\-\Up_*B_\r)\e
where $C\cap B_\r\-B_{\ep R}$ is identified with the graph of $d\xi$ for some smooth function $\xi:C\cap B_\r\-B_{\ep R}\to\R$ such that $\xi=\xi_0$ on $C\cap B_\r\-B_{2\ep R}$ and such that the graph of $d\xi$ is joined smoothly to $\ep\Ps_\ep(Y\cap K).$
%
%with $|f_\ep|_{C^k}=O(\ep^\la),$ $k\in\Z_{\ge0},$ on $B_R$
%If $k\in\Z_{\ge0}$ then $|f_\ep|_{C^k}=O(\ep)$ relative to the cylindrical metric on $C'.$
%This $f_\ep$ extends smoothly to $C\cap B_\r\-B_{\ep R}$ with $f_\ep=\xi_0$ on $C\cap B_\r\-B_{2\ep R}$ and with
%$|f_\ep|_{C^k}=O(\ep^\la)$ for every $k\in\Z_{\ge0}.$
%

Clearly, $Z$ is singular only at one point $p_Y\in Y\cap K;$ and we put $Z':=Z\-\{p_Y\}.$
Since $\Ps_\ep$ is at $p_Y\in\CP^1$ the identity up to the linear order it follows that $\Ps_\ep^*\ep^*\Om_1=\ep^2\Th_1$ at $p_Y.$
The tangent cone to $Z$ at $p_Y$ is accordingly a special Lagrangian cone with respect to $(\om_1,\Om_1)$ and isomorphic to Harvey--Lawson's $T^2$-cone.

\subsection{Perturbation to the True Solution}\l{13}
We formulate and prove Theorem \ref{S1} by which we perturb $Z$ to the special Lagrangian $X_1$ in $(N_1;\om_1,\Om_1).$
We use therefore the method by Joyce \cite{J1,J2,J3} but together with the extension by Pacini \cite{P3} because $Z$ is singular.

We write $x\les y$ if $x,y\in\R$ with $x\le cy$ for some $c>0$ independent of $\ep,$ and $x\px y$ if both $x\les y$ and $y\les x.$

We introduce conical and cylindrical metrics on $X_0':=X_0\-\{n_0\},$ $Y'$ and $Z',$ respectively. 
The {\it conical} metrics on $X_0',$ $Y'$ and $Z'$ shall mean those metrics induced from the K\"ahlerian ones on $(N',\om_0),$ $( Q_1,\th_1),$ and $( N_1',\om_1),$ respectively.
Take an $\ep$-independent smooth function $r_0:X_0'\to\R_+$ which agrees near $n_0$ with $r_C$ under $C'\cap B_\r\cong\gr d\xi_0.$
Take an $\ep$-independent smooth function $r_Y:Y'\to\R_+$ which agrees on $Y\-K$ with $\pi^*r_C,$ and near $p_Y$ with the distance from $p_Y$ relative to the conical metric of $Y'.$
Take a smooth function $r:Z'\to\R_+$ which satisfies, under \eq{Z}, $r\=r_0$ on $X_0\-\Up_*B_\r;$ $r\=r_C$ on $C\cap B_\r\-B_{\ep R};$ and, for each $k\in\Z_{\ge0},$ $|r|_{C^k}\px\ep r_Y,$ uniformly on $Y'\cap K$ with respect to its cylindrical metric.
The {\it cylindrical} metrics on $X_0',$ $Y'$ and $Z'$ shall mean the conformal transforms of their conical metrics by $r_0^{-2},$ $r_Y^{-2}$ and $r^{-2},$ respectively.

We introduce weighted H\"older spaces over $Z'.$
Take an $\ep$-independent manner of assigning unweighted H\"older norms of any functions on Riemannian manifolds.
Take $h\in(0,1)$ and $\mu\in(2,3);$ the latter we suppose so close to 2 as Joyce \cite[Definition 3.6]{J1} does, applying his definition to Harvey--Lawson's $T^2$-cone.
We denote by $C^{2,h}_\mu(Z';\R)$ the space of those functions $u:Z'\to\R$ with
$\|r^{-\mu}u\|_{C^{2,h}(Z')}<\iy$
where the H\"older norm respects the cylindrical metric on $Z'.$
We denote by $C^{0,h}_{\mu-2}(Z';\R)$ the same but with $(0,h;\mu-2)$ in place of $(2,h;\mu).$

We choose particular norms on $C^{2,h}_\mu(Z';\R)$ and $C^{0,h}_{\mu-2}(Z';\R).$
Take $\nu\in(2-\la,0).$
Take a smooth function $r^{\mu,\nu}:Z'\to\R_+$ which agrees, under \eq{Z}, with $r_0^\nu$ on $X\-\Up_*B_\r,$ and with $r_C^\nu$ on $C\cap B_\r\-B_{\ep R};$ and such that $|r^{\mu,\nu}|_{C^k}\px\ep^\nu r_Y^\mu,$ $k\in\Z_{\ge0},$ with respect to the cylindrical metric on $Y'\cap K,$ uniformly at every point of $Y'\cap K.$
We give $C^{2,h}_\mu(Z';\R)$ the norm $\|\bu\|_\nu$ by setting, for every $u\in C^{2,h}_\mu(Z';\R),$
$\|u\|_\nu:=\|(r^{\mu,\nu})^{-1}u\|_{C^{2,h}(Z')}.$
In the same manner but with $(0,h;\mu-2,\nu-2)$ in place of $(2,h;\mu,\nu),$ we take $r^{\mu-2,\nu-2}:Z'\to\R_+$ and give $C^{0,h}_{\mu-2}(Z';\R)$ the norm $\|\bu\|_{\nu-2}.$

We introduce Weinstein neighbourhoods of $X_0',$ $Y'$ and $Z'$ in $(N',\om_0),$ $( Q_1,\th_1),$ and $( N_1,\om_1),$ respectively.
There exists, as Joyce \cite[Theorem 4.6]{J1} constructs, a Weinstein neighbourhood of $X_0'$ in $(N',\om_0)$ which is, near $n_0,$ approximately dilation-invariant.
There exists likewise a Weinstein neighbourhood of $Y'$ in $( Q_1,\th_1)$ which is, near $p_Y,$ approximately dilation-invariant and is, near infinity, contained in the dilation-invariant Weinstein neighbourhood of $C'$ in $(Q_0',\th_0).$
From these two Weinstein neighbourhoods one gets, as Joyce \cite[Definition 6.7]{J3} does, a Weinstein neighbourhood of $Z'$ in $( N_1,\om_1)$ such that if
\[u\in\cU:=\{u\in C^{2,h}_\mu(Z';\R):\|u\|_\nu\le\ep^\la\}\]
then the graph of $du$ lies in the Weinstein neighbourhood of $Z'.$
%Given $u\in C^\iy(Z';\R)$ with $\gr du$ contained in this Weinstein neighbourhood of $Z'$ then $\gr du$ over $\ep\Ps_\ep(Y'\cap K)$ corresponds to the graph of $\ep^{-2}\Ps_\ep^*\ep^*du$ over $Y'\cap K;$
%the latter we abbreviate into $\ep^{-2}du$ by an abuse of notation.

%The graph of $du$ over $X_0\-\Up_*B_\r\sb Z'$ lies in the Weinstein neighbourhood of $Z'$ because $|du|\les\ep^\la$ where $|\bu|$ respects the cylindrical metric on $X_0\-\Up_*B_\r.$
%Turning to $\gr df_\ep\cong C\cap B_\r\-B_{\ep R}$ we have
%$r_C^{-2}|du|\le r_C^{\nu-2}\|u\|_\nu\les\ep^{\la+\nu-2}$
%where $|\bu|$ respects the cylindrical metric on $C\cap B_\r\-B_{\ep R}.$
%If $s\in[\ep R,\r]$ then $\gr du$ over $Z'\cap\d B_r$ corresponds to the graph of  $d(u+f_\ep)$ over $C\cap\d B_s.$
%This graph corresponds in turn to the graph of $s^{-2}\dl_s^*d(u+f_\ep)$ over $C\cap\d B_1$ which graph, because $r_C^{-2}|du|\les\ep^{\la+\nu-2},$ lies in the dilation-invariant Weinstein neighbourhood of $C'.$
%The graph of $du$ over $\gr df_\ep\sb Z'$ lies thus in the Weinstein neighbourhood of $Z'.$
%Finally, the graph of $\ep^{-2}du$ over $Y'\cap K\sb Z'$ lies in the Weinstein neighbourhood of $Z'$ because
%$\ep^{-2}|du|\le\ep^{-2}\ep^\nu r_Y^\mu\|u\|_\nu\les\ep^{\la+\nu-2}r_Y^\mu$
%where $|\bu|$ respects the cylindrical metric on $Y'\cap K;$ q.e.d.

We introduce, using the method by Joyce \cite[\S5.1]{J2}, some perturbations of the tangent cone to $Z$ at $p_Y.$
There are a unique smooth function $\al_1: N_1'\to\R_+$ such that $\om_1^3=(3i/4)\al_1^2\Om_1\w\ov\Om_1,$ and a principal $SU_3$-bundle over $K\sb Q_1$ whose fibre over a point $p$ consists of those linear isomorphisms of $T_{p_Y} Q_1$ onto $T_p Q_1$ which preserve both $\th_1,\al_1\Om_1$ evaluated at the two points $p_Y,p$ respectively.
%Take a sufficiently small slice of this bundle at the origin, which corresponds to the identity-map of $T_{p_Y} Q_1,$ relative to the maximal torus $T^2\sb SU_3$ (which preserves Harvey--Lawson's $T^2$-cone).
%This slice is diffeomorphic to a neighbourhood of 0 in the Lie algebra $\C^3\op\fsu_3/\ft^2$ which we denote by $V.$
%Recalling that every affine transformation of $\C^3$ is a Hamiltonian diffeomorphism, cutting off those Hamiltonians corresponding to the elements of $V,$ and making $V$ smaller if need be,
There are moreover a neighbourhood $V$ of 0 in $\C^3\op\fsu_3/\ft^2,$ and a smooth family $(\Ga_v)_{v\in V}$ of Hamiltonian diffeomorphisms of $ Q_1$ with Hamiltonians supported near $p_Y$ and with $d\Ga_v|_{p_Y}=v$ as an element of the $SU_3$-bundle.
We denote by $\spt\Ga_v,$ $v\in V,$ the compact supports of the Hamiltonians.

Thus, if $(u,v)\in\cU\t V$ then $\Ga_v(\gr du)$ is a compact Lagrangian in $( N_1,\om_1)$ with one-point singularity at $\Ga_v(p_Y)$ modelled upon Harvey--Lawson's $T^2$-cone.
This $\Ga_v(\gr du)$ is, outside its singular point, though a-priori only $C^{1,h}$-smooth yet, if special Lagrangian, $C^\iy$-smooth; and we prove indeed:
\begin{thm}\l{S1}
There exists $(u,v)\in\cU\t V$ such that $X_1:=\Ga_v(\gr du)$ is special Lagrangian with respect to $\Om_1.$
\end{thm}
We put, for every $(u,v)\in\cU\t V,$
\[\cP(u,v):=* (du)^*\Ga_v^*\Im\al_1\Om_1\]
where the first $*$ denotes Hodge's operator over $Z'$ with respect to its conical metric.
%where $du$ maps $Z'$ into its Weinstein neighbourhood.
Theorem \ref{S1} is then equivalent to:
\begin{thm}\l{S2}
There exists $(u,v)\in\cU\t V$ such that $\cP(u,v)=0.$
\end{thm}

The rest of this section is devoted to the proof of Theorem \ref{S2}.

Some results by Joyce \cite[\S6]{J2}, who deals with compact special Lagrangians with isolated conical singularities, extend immediately to compact nearly-special Lagrangians with isolated conical singularitites.
Analogy to the result by Joyce \cite[Proposition 6.4]{J2} implies that $\cP$ is a smooth map of $\cU\t V$ into $C^{0,h}_{\mu-2}(Z';\R).$

\begin{lem*}
The image of $\cP$ lies in
\[\cW:=\bigl\{w\in C^{0,h}_{\mu-2}(Z';\R):\ts\int_{Z'}\al_1|_{Z'}w*1=0\bigr\}.\]
\end{lem*}
\begin{proof}
The definitions of $\cP$ and $\al_1$ imply that if $(u,v)\in\cU\t V$ then
\e\l{z1}\ts\int_{Z'}\al_1|_{Z'}*\cP(u,v)=\Ga_{v*}(du)_*[Z]\cdot[\Im\Om_1]=[Z]\cdot[\Im\Om_1]\e
where $[Z]\in H_3( N_1';\R)$ and $[\Im\Om_1]\in H^3( N_1';\R).$
% where $ N_1'$ is used only in order that $H^3( N_1';\R)$ be a de Rham group.
Recall that $[Z]$ is independent of $\ep$ and that $[Z]\cdot[\Im\Om_1]$ converges as $\ep\to+0$ to $[X_0]\cdot[\Im\Om_0],$ which vanishes because $X_0$ is a special Lagrangian in $(N;\om_0,\Om_0).$
One finds then from \eq{z1} that $\cP(u,v)\in\cW.$
\end{proof}

Since Harvey--Lawson's $T^2$-cone is stable in the terminology by Joyce \cite[Definition 3.6 and \S3.2]{J2}
%\cite[Definition 3.4 and Example 3.5]{J5}
one gets, as Joyce \cite[\S6.2]{J2} does, a topological linear isomorphism
\[\cL:C^{2,h}_\mu(Z';\R)\op(\C^3\op\fsu_3/\ft^2)\to\cW,\]
which is the differential of $\cP$ at $0\in\cU\t V.$

We study now how $\cP$ depends upon $\ep.$
The main results are Corollary \ref{C1}, Theorem \ref{C2}, and Theorem \ref{UB}; which estimate $\cP(0),$ estimate $\cP-\cP(0)-\cL,$ and prove that $\cL$ is $\ep$-uniformly invertible, respectively.
These three resemble those by Pacini \cite[Proposition 6.2 (2)]{P3}, \cite[Proposition 5.9 (4)]{P3} and \cite[Theorem 4.8]{P3}, respectively.
Our circumstances are, however, more complicated in the respect that our ambient spaces are not flat whereas Pacini's are; we use therefore the auxiliary estimates in Lemma \ref{Xi}.

\begin{lem}\l{Xi}
Let $k,m\in\Z_{\ge0};$ and $E=E_\ep$ a smooth section of $\bot^mT^* Q_1$ over $K$ which vanishes at $p_Y$ and depends smoothly upon $\ep$ near and at $0.$
Then:
\item
{\bf(a)} $|r_Y^{-m}E|_{Y'}|_{C^{k,h}}\les r_Y;$
\item
{\bf(b)} if moreover $E_0=0$ then
$|r_Y^{-m}E|_Y'|_{C^{k,h}}\les\ep r_Y;$ and
\item
{\bf(c)} if $u\in\cU$ then
$|r_Y^{-m}(\ep^{-2}du)^*E|_{C^{0,h}}\les r_Y$
uniformly with respect to $u;$
\item
where all the estimates hold pointwise on $Y'\cap K,$ all the H\"older norms respecting the cylindrical metric on $Y'.$
\end{lem}
\begin{proof}
It is straightforward, which we leave to the reader; but we remark that (b) follows from (a) and that (c) extends (a) with $k=0.$
\end{proof}

\begin{thm*} One has, pointwise on $Z',$ for every $k\in\Z_{\ge0},$
\e\l{e0}
|\cP(0)|_{C^{k,h}}\les\begin{cases}
0\text{ wherever }r\ge2\ep R\\
\ep r_Y\text{ wherever }r\le\ep R\\
\ep^{\la-2}\text{ wherever }r\in[\ep R,2\ep R]
\end{cases}\e
with respect to the cylindrical metric on $Z'.$
\end{thm*}
\begin{proof}
The first line follows since that region defined by $r\ge2\ep R$ is a special Lagrangian in $(N;\om_0,\Om_0).$

We prove the second line of \eq{e0}.
Put $\Xi:=\ep^{-3}\Ps_\ep^*\ep^*\al_1\Om_1.$
There is a unique smooth function $\be_1: Q_1\to\R_+$ such that
$\th_1^3=(3i/4)\be_1^2\Th_1\w\ov\Th_1.$
Then, wherever $r\le\ep R$ or more generally on $Y'\cap K,$
\e\l{e1}
|\cP(0)|_{C^{k,h}}\px|r_Y^{-3}\Im\Xi|_{Y'}|_{C^{k,h}}=|r_Y^{-3}\Im(\Xi-\be_1\Th_1)|_{Y'}|_{C^{k,h}},
\e
the latter two norms respecting the cylindrical metric on $Y'\cap K.$
Applying Lemma \ref{Xi} (b) to $E:=\Im(\Xi-\be_1\Th_1)$ one gets from \eq{e1} the second line of \eq{e0}.

We prove the third line of \eq{e0}.
There is a unique smooth function $\be_0:Q_0'\to\R_+$ such that
$\th_0^3=(3i/4)\be_0^2\Th_0\w\ov\Th_0.$
One has then, on $C\cap B_{2R}\-B_R,$ putting $\ph:=\ep^{-2}\ep^*d\xi,$
\e\l{e3}
|\cP(0)|_{C^{k,h}}\px|\ph^*(\ep^{-3}\ep^*\Up^*\Im\be_0\Th_0)|_{C^{k,h}},
\e
the latter norm respecting the cylindrical metric on $C\cap B_{2R}\-B_R.$

We prove, on $C\cap B_{2R}\-B_R,$
\e\l{e4}
|\ph^*(\ep^{-3}\ep^*\Up^*\Im\be_0\Th_0-\Im\be_0\Th_0)|_{C^{k,h}}\les\ep.
\e
Using the method by Chan \cite[Theorem 4.9]{Ch1} and recalling $|\Up(x)-x|=O(|x|^2)$ one finds
\e\l{e6}
|\ep^{-3}\ep^*\Up^*\Im\be_0\Th_0-\Im\be_0\Th_0|_{C^{k,h}}\les\ep
\e
on $B_{2R}\-B_R$ with respect to its cylindrical metric.
On the other hand, since $|\xi|_{C^{k+2,h}}\les\ep^\la$ it follows that $|\ph|_{C^{k+1,h}}\les\ep^{\la-2}\les1$ on $C\cap B_{2R}\-B_R.$ This combined with \eq{e6} implies \eq{e4}.

We prove, still on $C\cap B_{2R}\-B_R,$
\e\l{e42}|\ph^*\Im \be_0\Th_0|_{C^{k,h}}\les\ep^{\la-2}.\e
We put, for every 1-form $\ps$ on $C\cap B_{2R}\-B_R,$ $\cG(\ps):=\ps^*\Im \be_0\Th_0.$
Since $C'$ is a special Lagrangian in $(Q_0';\th_0,\Th_0)$ it follows that $\cG(0)=0.$
Hence, recalling that $\cG(\ps)$ is pointwise on $C\cap B_{2R}\-B_R$ a smooth function of $(\ps,\nb\ps),$ one finds $|\cG(\ps)|_{C^{k,h}}\les|\ps|_{C^{k,h}}+|\nb\ps|_{C^{k,h}}.$
This applied to $\ps=\ph$ and combined with $|\ph|_{C^{k+1,h}}\les\ep^{\la-2}$ implies \eq{e42}.

The three estimates \eq{e1}, \eq{e4} and \eq{e42} imply the third line of \eq{e0}.
\end{proof}
\begin{cor}\l{C1}
$\|\cP(0)\|_{\nu-2}\les\ep^{\la-\nu}.$
\end{cor}
\begin{proof}
The definition of $\|\bu\|_{\nu-2}$ and the first line of \eq{e0} imply
\[\|\cP(0)\|_{\nu-2}\px
\ep^{2-\nu}\|r_Y^{2-\mu}\cP(0)\|_{C^{0,h}(Y'\cap K)}
+\ep^{2-\nu}\|\cP(0)\|_{C^{0,h}(C\cap B_{2\ep R}\-B_{\ep R})}
\]
with respect to the cylindrical metrics on $Y'\cap K$ and $C\cap B_{2\ep R}\-B_{\ep R}$ respectively.
The second line of \eq{e0} implies $\|r_Y^{2-\mu}\cP(0)\|_{C^{0,h}(Y'\cap K)}\les\|r_Y^{3-\mu}\|_{C^{0,h}(Y'\cap K)}\les1\les\ep^{\la-2}$ because $\mu<3,$ $\la<3.$
The third line of \eq{e0} implies $\|\cP(0)\|_{C^{0,h}(C\cap B_{2\ep R}\-B_{\ep R})}\les\ep^{\la-2}.$
These estimates imply Corollary \ref{C1}.
\end{proof}

We put $\cQ:=\cP-\cP(0)-\cL;$ and, taking an $\ep$-independent norm on $\C^3\op\fsu_3/\ft^2,$ for every $u\op v\in C^{2,h}_\mu(Z';\R)\op(\C^3\op\fsu_3/\ft^2),$
\[\|u\op v\|_\nu:=\|u\|_\nu+\ep^{2-\nu}|v|.\]
\begin{thm}\l{C2}
$\|\cQ z_1-\cQ z_0\|_{\nu-2}\les\ep^{\nu-2}\|z_1-z_0\|_\nu(\|z_1\|_\nu+\|z_0\|_\nu)$ uniformly with respect to $z_0,z_1\in\cU\t V;$ which uniformity we include in $\les,\px$ below.
\end{thm}
\begin{proof}
We prove that if $u_0,u_1\in\cU$ and $v\in V$ then
\e\l{Q1}
|\cQ(u_1,v)-\cQ(u_0,v)|_{C^{0,h}}\les r^{-4}|u_1-u_0|_{C^{2,h}}(|u_1|_{C^{2,h}}+|u_0|_{C^{2,h}})
\e
uniformly at every point of $Z,$ all the H\"older norms respecting the cylindrical metric on $Z.$ We prove moreover that if $u\in\cU$ and $v_0,v_1\in V$ then $\cQ(u,v_1)-\cQ(u,v_0)$ is supported in $\spt\Ga_v$ with
\e\l{Q2}|\cQ(u,v_1)-\cQ(u,v_0)|_{C^{0,h}}\les r_Y|v_1-v_0|(|v_1|+|v_0|),\e
the H\"older norm respecting again the cylindrical metric on $Z.$

Letting $p\in Z'$ and $\ul{ u}:=du|_p\op\nb du|_p$ we treat $\cP(u,v),$ $\cQ(u,v)$ as smooth functions of $\ul{ u}\op v\in[T_p^*Z'\op(T_p^*Z')^{\ot2}]\op(\C^3\op\fsu_3/\ft^2).$
One may then write
$\cQ(\ul{ u}\op v)=\cP(\ul{ u}\op v)-\cP(0)-\d\cP(0)(\ul{ u}\op v).$
We prove
\e\l{-4}|\d_{\ul{ u}}^2\cP|_{C^{0,h}}\les r^{-4}(p).\e

If $r(p)\ge\r$ then $\cP$ is independent of $\ep;$ which implies $|\d_{\ul{ u}}^2\cP|_{C^{0,h}}\les1$ and, accordingly, \eq{-4}.

Suppose next $r(p)\in[\ep R,\r].$ Let $s>0$ with $r(p)\in[s,2s]\sb[\ep R,\r].$
Put
\[\cG_s(\ps):=\ps^*(s^{-3}s^*\Up^*\Im\be_0\Th_0)\]
where $\ps$ is a 1-form on $C\cap B_2\-B_1.$
Regard $\cG_s(\ps)$ as a function of $\ul{\ps}:=\ps\op\nb\ps.$
One has then, with respect to the cylindrical metric on $C'\cap B_2\-B_1,$
\e\l{A1}
|\d_{\ul{ u}}^2\cP|_{C^{0,h}}\px|\d_{\ul{ u}}^2\cG_s(s^{-2}s^*du)|_{C^{0,h}}
\px
s^{-4}|\d_{\ul{\ps}}^2\cG_s(s^{-2}s^*du)|_{C^{0,h}}
\e
Since $s^{-3}s^*\Up^*\Im\be_0\Th_0$ converges smoothly as $s\to+0$ over $B_2\-B_1$ it follows that $|\d_{\ul{\ps}}^2\cG_s|_{C^{0,h}}\les1.$ This combined with \eq{A1} implies
$|\d_{\ul{ u}}^2\cP|_{C^{0,h}}\les s^{-4}\px r^{-4}(p)$ and, accordingly, \eq{-4}.

Suppose next $r(p)\le\ep R$ with $p\notin\spt\Ga_v.$
Put
\[\cF_\ep(\ps):=\ps^*\Im\al_1\Om_1\]
where $\ps$ is a 1-form on $Y'\cap K.$
Regard $\cF_\ep(\ps)$ as a function of $\ul{\ps}:=\ps\op\nb\ps.$
One has then, with respect to the cylindrical metric on $Y'\cap K,$
\e
\l{A2}
|\d_{\ul{ u}}^2\cP|_{C^{0,h}}\px|\d_{\ul{ u}}^2\cF_\ep(\ep^{-2}du)|_{C^{0,h}}
\px
\ep^{-4}|\d_{\ul{\ps}}^2\cF_\ep(\ep^{-2}du)|_{C^{0,h}}.
\e
Since $\cF_\ep$ is smooth at $\ep=0$ it follows that
$|\d_{\ul{\ps}}^2\cF_\ep|_{C^{0,h}}\les1.$
This combined with by \eq{A2} implies $|\d_{\ul{ u}}^2\cP|_{C^{0,h}}\les\ep^{-4}\px r^{-4}(p)$ and, accordingly, \eq{-4}.

Combining the techniques in the previous two cases, one gets \eq{-4} even if $p\in\spt\Ga_v.$

In general, let $F$ be a smooth function on a Banach space $X;$ $Q:=F-F(0)-dF(0);$ and $x_0, x_1\in X;$ then,
\e\l{QQ}Q(x_1)-Q(x_0)=(x_1-x_0)\lrcorner\int_0^1 x_t\lrcorner\int_0^1\d^2F(t' x_t)dt'dt\e
where $x_t:=x_0+t(x_1-x_0).$
This applied to $\cP$ over $p$ and combined with $|\d_{\ul{ u}}^2\cP|_{C^{0,h}}\les r^{-4}(p)$ implies \eq{Q1}.

As to \eq{Q2}, the definition of $\Ga_v$ implies $\cQ(u,v_1)-\cQ(u,v_0)$ supported in $\spt\Ga_v.$
One has, putting $E:=\d_v^2\Ga_v^*\Im\al_1\Om_1,$
\e\l{QQ1}|\d_v^2\cP|_{C^{0,h}}\px|r_Y^{-3}(\ep^{-2}du)^*E|_{C^{0,h}}\e
where the latter norm respects the cylindrical metric on $Y'.$
Applying Lemma \ref{Xi} (c) to $E$ one finds from \eq{QQ1} that
$|\d_v^2\cP|_{C^{0,h}}\les r_Y.$
This combined with \eq{QQ} implies \eq{Q2}.

The two estimates \eq{Q1}, \eq{Q2} and the definitions of $\|\bu\|_\nu,\|\bu\|_{\nu-2}$ imply Theorem \ref{C2}.
\end{proof}

\begin{thm}\l{UB}
$\cL:C^{2,h}_\mu(Z';\R)\op(\C^3\op\fsu_3/\ft^2)\to\cW$
is $\ep$-uniformly invertible.
\end{thm}
\noindent
The proof of this theorem, given shortly, is rather different from that by Pacini \cite[Proof of Theorem 4.8]{P3}.
The auxiliary factor $\ze$ in our proof is unnecessary to Pacini's; for in his circumstances, the relevant linear operators are all invertible because of the noncompactness of special Lagrangians.
We use therefore the more general method by Donaldson and Kronheimer \cite[\S7.2.2]{DK} in Yang--Mills gauge theory.
\begin{proof}[Proof of Theorem $\ref{UB}$]
Take $\ze\in C^\iy(Z';\R)$ supported in $X\-\Up_*B_\r$ with $\int_{Z}\ze*1=1.$
The definition of $\cW$ implies that the operator
\[\ze\op\cL:\R\op[C^{2,h}_\mu(Z';\R)\op(\C^3\op\fsu_3/\ft^2)]\to C^{0,h}_{\mu-2}(Z';\R)\]
defined by $t\op z\mapsto t\ze+\cL z$ is, with $\ep$ fixed, a topological linear isomorphism.
It suffices then to prove that $\ze\op\cL$ is $\ep$-uniformly invertible.

There is a unique smooth function $t:Z'\to\R/2\pi\Z$ with $\al_1\Om_1|_{Z'}=e^{it}*1.$
Standard computation implies $\cL u=-d^*(\cos t\,du),$ $u\in C^{2,h}_\mu(Z';\R),$ where $d^*$ respects the conical metric on $Z'.$

We approximate $\cL$ by two $\ep$-independent operators over $X_0',$ $Y'$ which we define shortly in \eq{LX}, \eq{LY} respectively.
We construct from the latter two operators an approximate right inverse to $\ze\op\cL,$ which is $\cK_0+\cK_Y$ in our notation.
The definitions of $\cK_0,$ $\cK_Y$ use some cut-off functions.
Lemmata \ref{d1} and \ref{d2} prove the relevant estimates for $\cK_0$ and $\cK_Y$ respectively.

%\frac{d}{ds}(s\,du)^*\Im \Xi_\ep|_{s=0}=d(J\grad_\ep u\lrcorner\Im \Xi_\ep|_{Z'})
%=d(\grad_\ep u\lrcorner \Im i\Xi_\ep|_{Z'})
%=d(\grad_\ep u\lrcorner \cos t*1)=
%-d^*(\cos t\,du)$
%where $J$ denotes the complex structure of $ N_1'$ and where $\grad_\ep,d_\ep^*$ respect the conical metric on $Z'.$
We denote by $C^{2,h}_\nu(X_0';\R)$ the Banach space of those $u\in C^{2,h}(X_0';\R)$ with
\[\|u\|_{C^{2,h}_\nu(X_0')}:=\|r_0^{-\nu}u\|_{C^{2,h}(X_0')}<\iy\]
relative to the cylindrical metric on $X_0'.$
We denote by $C^{0,h}_{\nu-2}(X_0';\R)$ the same but with $(0,h;\nu-2)$ in place of $(2,h;\nu).$

Define
\e\l{LX}\ze\op\cL_0:\R\op C^{2,h}_\nu(X_0';\R)\to C^{0,h}_{\nu-2}(X_0';\R)\e
by $t\op u\mapsto t\ze-d_0^*du$ where $d_0^*$ respects the conical metric on $X_0'.$
Joyce \cite[Theorems 2.14 and 2.16 (b)]{J1} proves that $\ze\op\cL_0$ is surjective; and accordingly, admits a bounded right inverse $\pi_0\op\cR_0.$

Take an $\ep$-independent smooth function $r_Y^{\mu,\nu}:Y'\to\R_+$  with $r_Y^{\mu,\nu}\=r_Y^\nu$ on $Y\-K$ and with $r_Y^{\mu,\nu}\=r_Y^\mu$ near $p_Y.$
We denote by $C^{2,h}_{\mu,\nu}(Y';\R)$ the Banach space of those $u\in C^{2,h}(Y';\R)$ with
\[\|u\|_{C^{2,h}_{\mu,\nu}(Y')}:=\|(r_Y^{\mu,\nu})^{-1}u\|_{C^{2,h}(Y')}<\iy\]
relative to the cylindrical metric on $Y'.$
We denote by $C^{0,h}_{\mu-2,\nu-2}(Y';\R)$ the same but with $(0,h;\mu-2,\nu-2)$ in place of $(2,h;\mu,\nu).$

Take an $\ep$-independent $\et\in C^\iy(Y';\R)$ which vanishes near $p_Y$ and is constant outside $K.$
Define
\e\l{LY}\cL_Y:C^{2,h}_{\mu,\nu}(Y';\R)\op\R\op(\C^3\op\fsu_3/\ft^2)\to C^{0,h}_{\mu-2,\nu-2}(Y';\R)\e
by $u\op t\op v\mapsto-d_Y^*d(u+t\et)+\ga _0v$ where $d_Y^*$ respects the conical metric on $Y';$ and $\ga_0$ is derived from $v\mapsto*_Y(\Ga_v|_{\ep=0})^*\Im\be_1\Th_1$ at $v=0,$ with $*_Y$ respecting the conical metric on $Y'.$
Pacini \cite[Theorem 6.10]{P3} proves that $\cL_Y$ is a topological linear isomorphism.

We introduce now the cut-off functions.
There are, with $K$ and $R$ made larger enough if need be, two $\ep$-independent constants $S,T>0$ such that $2Se^T<R$ and \eq{Z} makes sense with $e^{-T}S$ in place of $R.$
Let $F_0,F_Y$ be a partition of unity on $Z'$ subordinate to its open covering by $r>\ep S,$ $r<2\ep S$ respectively.
Put $T:=\log(S/R).$
Let $\ep$ be so small that $2\ep S e^T<\r.$
Take $G_0,G_Y\in C^\iy(Z';\R)$ with
\[G_0\=\begin{cases}1\text{ wherever }r\ge\ep e^{-T}S\\
0\text{ wherever }r\le\ep S,
\end{cases}
G_Y\=\begin{cases}1\text{ wherever }r\le 2\ep S\\
0\text{ wherever }r\ge\ep R,
\end{cases}\]
and $|G_0|_{C^{2,h}}+|G_Y|_{C^{2,h}}\les T^{-1}$ relative to the cylindrical metric on $Z'.$

Define
\[\cK_0:C_{\mu-2}^{0,h}(Z';\R)\to \R\op C^{2,h}_\mu(Z';\R).\]
by $\cK_0:=(1\op G_0)(\pi_0\op\cR_0)F_0.$
Define
\[\cK_Y:C_{\nu-2}^{0,h}(Z';\R)\to C^{2,h}_\mu(Z';\R)\op(\C^3\op\fsu_3/\ft^2)\]
by $\cK_Y:=G_Y\cL_Y^{-1}F_Y$ where $G_Y$ acts linearly upon $C^{2,h}_{\mu,\nu}(Y';\R)\op\R$ as $u\op t\mapsto G_Y(u+t\et)\in C^{2,h}_\mu(Z';\R)$ and upon $\C^3\op\fsu_3/\ft^2$ as the identity.
Then:
\begin{lem}\l{d1}
$\|(\ze\op\cL)\cK_0w-F_0w\|_{\nu-2}\les(T^{-1}+\ep^{\la-2})\|w\|_{\nu-2}$
uniformly with respect to $w\in C^{0,h}_{\mu-2}(Z';\R);$ which uniformity we include in $\les,\px$ below.
\end{lem}
\begin{proof}
Put $u:=\cR_0F_0w.$ Straightforward computation implies
\e\l{L0}(\ze\op\cL)\cK_0w-F_0w=[\cL_0,G_0]u+(\cL-\cL_0)G_0u.\e
We estimate $[\cL_0,G_0]u.$
There are over $\spt G_0$ some sections $C_1,C_2$ of $TX_0',TX_0'^{\ot2}$ respectively with $|C_1|_{C^{0,h}}+|C_2|_{C^{0,h}}\les1$ and
\e\l{d00}r_0^2[\cL_0,G_0]u=-\langle dG_0,du\rangle+u(C_1\nb G_0+C_2\nb^2G_0)\e
relative to the cylindrical metric on $X_0'.$
Since $|G_0|_{C^{2,h}}\les T^{-1}$ it follows, by \eq{d00}, that
\[
\|[\cL_0,G_0]u\|_{C_{\nu-2}^{0,h}(X_0')}\les T^{-1}\|u\|_{C_\nu^{2,h}(X_0')}.\]
Hence one finds, recalling that $u$ is supported in $\spt F_0,$
\e
\l{L1}\|[\cL_0,G_0]u\|_{\nu-2}\les T^{-1}\|w\|_{\nu-2}.
\e

We estimate $(\cL-\cL_0)G_0u,$ which is the last term on \eq{L0}.
Identify the two annular regions in $X$ and $Z$ of inner radius $\ep e^{-T}S$ and outer radius $\ep R$ with respect to $r_0,r$ respectively.
One has then, putting $\ph:=G_0u,$
\e\l{d02}\cL\ph=-d_\ep^*(\cos t\,d\ph)\text{ and }\cL_0\ph=-d_0^*d\ph.\e
Since $\cP(0)=\sin t$ it follows, by the third line of \eq{e0}, that
\e\l{d03}|\cos t-1|_{C^{1,h}}\les\ep^{\la-2}\e
on $\spt\ph$ with respect to the cylindrical metric on $X_0'.$
Let $g_0$ be this metric, and $g$ the cylindrical metric on $Z';$ then, analogy to the estimate by Joyce \cite[Theorem 5.2]{J1} implies, on $\spt\ph,$
\e\l{d04}|g-g_0|_{C^{1,h}}\les\ep^{\la-2}.\e
By \eq{d02}--\eq{d04}, there exist over $\spt\ph$ some sections $E_1,E_2$ of $TX_0',TX_0'^{\ot2}$ respectively with
\[(\cL-\cL_0)\ph=r_0^{-2}(E_1\nb\ph+E_2\nb^2\ph)\text{ and }|E_1|_{C^{0,h}}+|E_2|_{C^{0,h}}\les\ep^{\la-2}.\]
Consequently, $\|(\cL-\cL_0)G_0u\|_{\nu-2}\les\ep^{\la-2}\|w\|_{\nu-2};$ which combined with \eq{L0} and \eq{L1} implies Lemma \ref{d1}.
\end{proof}

\begin{lem}\l{d2}
$\|\cL\cK_Yw-F_Yw\|_{\nu-2}\les(T^{-1}+\ep)\|w\|_{\nu-2},$
$w\in C^{0,h}_{\mu-2}(Z';\R).$
\end{lem}
\begin{proof}
Put $y:=\cL_Y^{-1}F_Yw.$ The definition of $\cK_Y$ implies
\e\l{d20}
\cL\cK_Yw-F_Yw=[\cL_Y,G_Y]y+(\cL-\cL_Y)G_Yy.
\e
The method of the proof of \eq{L1} implies
\e\l{d21}
\|[\cL_Y,G_Y]y\|_{\nu-2}\les T^{-1}\|w\|_{\nu-2}.
\e
Put $y=:u\op t\op v,$ $\ph:=d[G_Y(u+t\et)]$ and $\ga:=\cL|_{\C^3\op\fsu_3/\ft^2}.$ Then,
\e\l{d31}
\cL G_Yy=-d^*(\cos t\,\ph)+\ga v\text{ and }
\cL_YG_Yy=-d_Y^*\ph+\ga_0v.
\e
The second line of \eq{e0} implies
\e\l{d312}|\cos t-1|_{C^{1,h}}\les\ep\e
on $Y'\cap K$ with respect to the cylindrical metric on $Y'.$
Letting $g_Y$ be this metric, letting $g$ be the cylindrical metric on $Z',$ and applying Lemma \ref{Xi} (b) to $r_Y^2(g_Y-g)$ on $Y'\cap K,$ one finds
\[|g_Y-g|_{C^{1,h}}\les\ep.\]
This combined with \eq{d312} implies
\e\l{d4}
|d^*(\cos t\,\ph)-d_Y^*\ph|_{C^{0,h}}\les\ep|\ph|_{C^{1,h}}.
\e

Let $\cH(v):=*\Ga_v^*\Im \Xi-*_Y(\Ga_v|_{\ep=0})^*\Im\be_1\Th_1;$ and $\cH'$ the differential of $\cH$ at $v=0;$ then, $\ga-\ga_0=\cH'|_{Y'}.$ Hence one finds, applying Lemma \ref{Xi} (b) to $\cH',$
\[|(\ga-\ga_0)v|_{C^{0,h}}\les\ep r_Y|v|.\]
This combined with \eq{d31} and \eq{d4} implies
$\|(\cL-\cL_Y)G_Yy\|_{\nu-2}\les\ep\|w\|_{\nu-2}.$
Consequently, by \eq{d20} and \eq{d21}, Lemma \ref{d2} holds.
\end{proof}

We complete the proof of Theorem $\ref{UB}.$
Put $\cJ:=(\ze\op\cL)(\cK_0+\cK_Y);$ then, Lemma \ref{d1}, Lemma \ref{d2}, and the identity $F_0+F_Y=1$ imply
\[\|\cJ w-w\|_{\nu-2}\les(\ep^{\la-2}+T^{-1})\|w\|_{\nu-2}.\]
Thus, with $T$ large enough and $\ep$ small enough, $\cJ$ is $\ep$-uniformly invertible; accordingly, $(\ze\op\cL)^{-1}=(\cK_0+\cK_Y)\cJ^{-1}$ is $\ep$-uniformly bounded; and consequently, $\cL^{-1}$ is also $\ep$-uniformly bounded.
\end{proof}

We complete the proof of Theorem $\ref{S2}.$
Solving the equation $\cP(u,v)=0$ is equivalent to finding a fixed point of
$\cT:=-\cL^{-1}[\cP(0)+\cQ].$
Let $\cB$ denote the closed ball of radius $\ep^\la$ with respect to $\|\bu\|_\nu$ in $C_\mu^{2,h}(Z';\R)\op(\C^3\op\fsu_3/\ft^2).$
The definition of $\cU\t V$ implies $\cB\sb\cU\t V.$
We prove that $\cT$ maps $\cB$ into itself.
Corollary \ref{C1}, Theorem \ref{C2} and Theorem \ref{UB} imply that if $x\in\cB$ then
\[\|\cT x\|_\nu\les(\ep^{-\nu}+\ep^{\la+\nu-2})\ep^\la.\]
Hence one finds, recalling $-\nu>0$ and $\la+\nu-2>0,$ and making $\ep$ small enough, $\cT x\in\cB.$
Theorems \ref{C2} and \ref{UB} imply moreover that if $x,y\in\cB$ then
\[\|\cT x-\cT y\|_\nu\les\ep^{\la+\nu-2}\| x-y\|_{\nu-2}.\]
Thus, with $\ep$ small enough, $\cT$ contracts $\cB;$ and consequently, there exists a unique $(u,v)\in\cB$ with $\cP(u,v)=0;$ q.e.d.

\section{Example}\l{ex}
This section contains four subsections: \S\ref{21} presents the example of the building blocks; \S\ref{22} recalls the relevant general facts from the preceding papers \cite{J5,I}; and \S\ref{23}, \S\ref{24} deal with the projective small resolutions and the smoothings, respectively, of the nodal complex $3$-fold $N_1$ in the example.

\subsection{Example of the Building Blocks}\l{21}
We construct in $\CP^4$ a quintic hypersurface $N_0.$
We denote by $(z_0,z_1,z_2,z_3,z_4)$ the homogeneous co\"ordinates of $\CP^4.$
Straightforward computation implies:
\begin{lem}\l{EN}{\bf(a)}
There are quintics of the form $z_1g_1-z_2g_2,$ where $g_1,g_2$ are both quartic homogeneous polynomials, which has at least sixteen singular points counted with multiplicities, defined by $z_1=z_2=g_1=g_2=0.$ With $g_1,g_2$ generic this quintic has exactly sixteen nodes.
\item
{\bf(b)}
There are quintics of the form $z_1z_4f_1-z_2z_3f_2,$ where $f_1,f_2$ are both homogeneous cubics, which has at least $49$ singular points counted with multiplicities, defined by two of $\{z_1,z_4,f_1\}$ and two of $\{z_2,z_3,f_2\}$ both zero.
One of the $49$ is defined by $(z_1,z_2,z_3,z_4)=0.$
On twelve of the $49,$ three of $\{z_1,z_2,z_3,z_4\}$ are zero; there are, according to the unique $z_i$ which does not vanish, four cases; in each of which, there are three singular points counted with multiplicities, defined by either $f_1=0$ or $f_2=0.$
On the remaining $36,$ one of $\{z_1,z_2\}$ and one of $\{z_3,z_4\}$ are zero; there are, according to those $z_i$ which do and do not vanish, four cases; in each of which, there are nine singular points counted with multiplicities, defined by $f_1=f_2=0.$
With $f_1,f_2$ generic this quintic has exactly $49$ nodes.
\end{lem}
We use not (a) but (b) of Lemma \ref{EN}.
Let $c\gg1,$ $f_1:=z_0^3-z_1^3-z_4^3,$ $f_2:=cz_0^3-z_2^3-z_3^3,$ and $N_0$ the quintic hypersurface in $\CP^4$ defined by $f_0:=z_1z_4f_1-z_2z_3f_2=0;$ then:
\begin{lem*}
$N_0$ is nodal with exactly $49$ nodes.
\end{lem*}
\begin{proof}
Let $n$ be a singular point of $N_0.$
We prove $z_1z_2z_3z_4=0$ at $n.$
One has, at $n,$
\begin{align*}
0=\d_1f_0\=z_4(z_0^3-4z_1^3-z_4^3)\text{ and }0=\d_4f_0\=z_1(z_0^3-z_1^3-4z_4^3).
\end{align*}
Accordingly, if $z_1z_2z_3z_4\ne0$ then $z_1^3=z_4^3=z_0^3/5;$ and in particular,
\e\l{5}|z_1z_4f_1|=(1/5)^{2/3}(3/5)|z_0|^5>0.\e
Likewise, since $\d_2f_0=\d_3f_0=0$ it would follow that $z_2^3=z_3^3=cz_0^3/5;$ and in particular,
$|z_2z_3f_2|=(c/5)^{2/3}(3c/5)|z_0|^5>0.$
This combined with \eq{5} contradicts $f_0=0$ with $c\ne1.$
Consequently, $z_1z_2z_3z_4=0$ at $n.$

We prove, in fact, $z_1z_4=0$ at $n.$ Otherwise, \eq{5} would hold again; and in particular, $z_1z_4f_1\ne0.$
This combined with the definition of $f_0$ implies $z_2z_3\ne0;$ which contradicts $z_1z_2z_3z_4=0.$
Consequently, $z_1z_4=0$ at $n.$

Likewise, $z_2z_3=0$ at $n.$

We prove that $n$ is a node.
One may suppose therefore, by symmetry, $z_1=z_2=0$ at $n.$
At $n,$ if $z_0=0$ then $\d_1f_0=\d_2f_0=0$ would imply $z_4=z_3=0,$ all the co\"ordinates vanishing; which is a contradiction and thus implies $z_0\ne0.$
Put $z_4':=z_4f_1$ and $z_3':=z_3f_2;$ then, $(z_1,z_2,z_3',z_4')$ with $z_0=1$ is a co\"ordinate-system about $n$ in $\CP^4.$ Since $f_0=z_1z_4'-z_2z_3'$ it follows that $n$ is a node.
This expression of $f_0$ implies moreover that every singular point of $N_0$ is defined as in Lemma \ref{EN} (b) so that $N_0$ has exactly 49 nodes.
\end{proof}

Using the inhomogeneous co\"ordinate-system $(z_1,z_2,z_3,z_4)$ with $z_0=1,$
let
\[\om_0:=(i/2)\d\db\log(1+|z_1|^2+c|z_2|^2+c|z_3|^2+|z_4|^2);\]
and $\Om_0$ the residue of the meromorphic 4-form $cf_0^{-1}dz_1\w dz_2\w dz_3\w dz_4.$
This $(\om_0,\Om_0)$ defines on $N_0$ an almost Calabi--Yau structure in the sense of Definition \ref{N}. 
Let $X_0$ be the fixed-point set in $N_0$ of the involution $(z_0,z_1,z_2,z_3,z_4)\mapsto(\bar z_0,\bar z_4,\bar z_3,\bar z_2,\bar z_1).$
This $X_0$ is certainly a special Lagrangian in $(N_0;\om_0,\Om_0)$ with one-point singularity at the node $n_0\in N_0$ defined by $(z_1,z_2,z_3,z_4)=(0,0,0,0).$
Let $D_0$ be the prime divisor on $N_0$ defined by $z_3=z_4=0.$
Then:
\begin{thm}\l{EB}
This $(N_0;\om_0,\Om_0;X_0;D_0)$ satisfies Hypothesis $\ref{H}.$
\end{thm}
\begin{proof}
Define a biholomorphism-of-germs $\Ph:(Q_0,0)\to(N,n_0)$ by $\Ph^{-1}=(x_1,x_2,x_3,x_4)$ where $x_1:=z_1,$ $x_2:=c^{1/2}z_2,$ $x_3:=c^{-1/2}z_3f_2,$ and $x_4:=z_4f_1,$ using the inhomogeneous co\"ordinate-system $(z_1,z_2,z_3,z_4)$ with $z_0=1.$
Then, $\Ph^*\om_0|_0=\th_0$ and, near $0$ in $Q_0,$ $\Ph^*\Om_0\=\Th_0.$
One has moreover, near 0 in $\Ph^*X_0',$ $\re x_2x_3>0;$ in which region, $\Ph^*X_0$ is defined by an equation of the form
\[(|x_1|^2-|x_4|^2,-|x_2|^2+|x_3|^2,\Im x_2x_3)=O(|x|^3).\]
Consequently, $\Ph^*X_0$ approaches $C$ with order $3.$
Finally, it is clear that $\Ph^*D_0\sb P_0.$ 
\end{proof}

\begin{lem}\l{hm}
{\bf(a)} $X_0'$ is homeomorphic to $S^1\t\R^2;$
and {\bf(b)} $H_1(X_0';\Z)$ is generated, under the natural homomorphism $H_1(C';\Z)\to H_1(X_0';\Z),$ by an $S^1$-orbit in $C'$ with respect to the action $(x_1,x_2,\bar x_2,\bar x_1)\mapsto(\ze x_1,x_2,\bar x_2,\ze^{-1}\bar x_1)$ where $\ze\in S^1$ and $(x_1,x_2,x_3,x_4)\in C'.$
\end{lem}
\begin{proof}
Every point of $X_0'$ is expressible as $(t,z_1,z_2,\bar{z}_2,\bar{z}_1)\in\CP^4$ with $t\in\R,$ $|z_1|^2+|z_2|^2=1,$ and
\[st^3=|z_1|^2(z_1^3+\bar{z}_1^3)-|z_2|^2(z_2^3+\bar{z}_2^3)\]
where $s:=|z_1|^2-c|z_2|^2.$ In particular,
\e\l{x12} |z_1|=\frac{c+s}{c+1}\text{ and }|z_2|=\frac{1-s}{1+c}.\e
We abbreviate $(t,z_1,z_2,\bar{z}_2,\bar{z}_1)$ to $(t,z_1,z_2),$ which expression is unique up to the relation $(t,z_1,z_2)=(-t,-z_1,-z_2).$
Take a sufficiently small positive constant $\si.$
Let $W_1,$ $W_2$ and $W_3$ be the three subsets of $X_0'$ defined by $s\le-\si,$ $|s|\le\si,$ and $s\ge\si,$ respectively.

We construct two homeomorphisms $W_1\to D^2\t S^1.$
Since $s\ne0$ on $W_1$ it follows that the projection $(t,z_1,z_2)\mapsto(z_1,z_2)$ embeds $X_1$ into $\RP^3.$
By \eq{x12} one has, on $W_1,$
\[|z_1|\le\frac{c-\si}{c+1}=:R_1>0\text{ and }|z_2|>0.\]
There is accordingly a homeomorphism
\e\l{ph0}(w,\ze):W_1\xrightarrow{\sim}D^2\t S^1\text{ where }w:=\frac{z_1z_2}{R_1|z_2|}\text{ and }\ze:=\frac{z_2^2}{|z_2|^2}.\e
Take a sufficiently small cut-off function $\ch:D^2\to[0,1]$ which is, near 0, a nonzero constant and, near $\d D^2,$ identically zero; then, there is another homeomorphism 
\e\l{ph1} W_1\xrightarrow{\sim}D_1^2\t S^1_1,\text{ defined by }(t,z_1,z_2)\mapsto\Bigl(\ze^{-1}w^2+\ch(w)w,\frac{z_1z_2}{|z_1z_2|}\Bigr),\e
where $D_1^2,$ $S^1_1$ are copies of $D^2,$ $S^1$ respectively.
We use \eq{ph1} rather than \eq{ph0}.

Put $A:=[-\si,\si]\t(\R\cup\{\iy\})$ and $A':=A\-\{(0,\iy)\}.$ We construct a homeomorphism $W_2\cong A'\t S^1_1.$
Take any point of $W_2;$ express it as $(t,z_1,z_2)$ with $\re z_1\ge0;$ and, wherever $s\ne0,$ let $m_+$ and $m_-$ denote the greater and less, respectively, of the two numbers
\[2s^{-1}[|z_1|^5-|z_2|^2\re(z_2^3)]\text{ and }-2s^{-1}[|z_1|^5+|z_2|^2\re(z_2^3)].\]
By \eq{x12}, with $\si$ small enough, $|z_1|$ is so close to 1 and $|z_2|$ is so close to 0 that $m_+>0,$ $m_-<0$ and $m_-\le t^3\le m_+.$
Define $\ta\in[-\iy,+\iy]$ by $\ta:=\tan(2t/\pi m_+)$ if $t\ge0,$ and $\ta:=\tan(2t/\pi m_-)$ if $t\le0.$
Identifying $\iy,-\iy$ one gets a map $[-\iy,+\iy]\to\R\cup\{\iy\}$ which we denote by $\ta\mapsto\ul\ta.$
Then, there is a homeomorphism
\[W_2\xrightarrow{\sim}A'\t S^1_1,\text{ defined by }(t,z_1,z_2)\mapsto\Bigl(s,\ul\ta;\frac{z_1z_2}{|z_1z_2|}\Bigr).\]
This combined with \eq{ph1} and a suitable homeomorphism $\d D^2_1\cong\{-\si\}\t(\R\cup\{\iy\})$ implies
\e\l{ph2}
W_1\cup W_2\cong (D^2_1\cup A')\t S^1_1\cong (D^2\-\{0\})\t S^1_1=:S^1_2\t(0,1]\t S^1_2=:S^1_2\t (D^2_2\-\{0\})\e
where $S^1_2,$ $D^2_2$ are other copies of $S^1,$ $D^2$ respectively.

By \eq{x12} one has, on $W_3,$
\[|z_1|>0\text{ and }|z_2|\le\frac{1-\si}{1+c}=:R_3>0.\]
There is accordingly a homeomorphism
\e\l{ph3}W_3\xrightarrow{\sim}S^1_2\t D^2_3,\text{ defined by }(t,z_1,z_2)\mapsto\Bigl(\frac{z_1^2}{|z_1^2|},\frac{z_1z_2}{|z_1|R_3}\Bigr),\e
where $D^2_3$ is yet another copy of $D^2.$
The two homeomorphisms \eq{ph2} and \eq{ph3} combined with another suitable homeomorphism $\d D_2^2\cong\d D_3^2$ imply
\[X_0'\cong S^1_2\t[(D^2_2\-\{0\})\cup D^2_3]\cong S^1\t\R^2\]
which completes the proof of (a).

As to (b), $H_1(X_0';\Z)$ is generated by a 1-cycle in $A'$ about its puncture with $z_2/|z_2|$ constant; on which, one may suppose $|t|=T,$ a sufficiently large constant. There exist $a_T^\pm,b_T^\pm\in\R$ with $a_T^\pm<0<b_T^\pm$ such that this 1-cycle is defined by $t=\pm T,$ $s\in[a_T^\pm,b_T^\pm]$ where $\pm$ is taken always consistently.
With $t=T,$ as $s$ moves from $b_T^+$ to $a_T^+$ the argument of $z_1$ moves from 0 to $\pi;$ whereas with $t=-T,$ as $s$ moves from $a_T^-$ to $b_T^-$ the argument of $z_1$ moves from 0 to $\pi.$
Hence one gets, recalling $(t,z_1,z_2)=(-t,-z_1,-z_2),$ that 1-cycle required in (b).
\end{proof}

\subsection{General Facts}\l{22}
Let $(M;\om,\Om)$ be an almost Calabi--Yau manifold of complex dimension 3, and $X$ a compact special Lagrangian in $M$ with one-point singularity $x$ modelled upon Harvey--Lawson's $T^2$-cone.

Suppose $X':=X\-\{x\}$ diffeomorphic to $S^1\t\R^2;$
which implies in particular, over $\Z,$ some natural isomorphisms
$H^1_c(X')=H^2(X')=0,$ $H^1(X')=H_c^2(X')=H^2(X,\{x\})=H^2(X)=\Z,$
and accordingly, a natural exact sequence
\e\l{ES}0\to H^1(X')\to H^1(T^2)\to H^2(X)\to0.\e

There are, up to dilations, exactly three models $\{L_1,L_2,L_3\}$ of the possible smoothings of $X;$ or more explicitly,
\[L_1:=\{(x_1,x_2,x_3)\in\C^3:|x_1|^2-1=|x_2|^2=|x_3|^2, x_1x_2x_3\in\R_{\ge0}\}\]
and $L_2,L_3$ are defined by cyclic permutations of $\{x_1,x_2,x_3\}$ in this expression.
These three are all noncompact closed smooth special Lagrangians in $\C^3$ diffeomorphic to $S^1\t\R^2$ and asymptotic at infinity to Harvey--Lawson's $T^2$-cone of multiplicity one.

There are also three vectors $y(L_i)\in H^1(T^2;\Z),$ $i\in\{1,2,3\},$ whose sum is zero; which are pairwise $\R$-linearly-independent; and which generate the images of the natural injections $H^1(L_i;\Z)\to H^1(T^2;\Z),$ respectively. In fact, $y(L_i)$ is proportional by a positive constant $($independent of $i)$ to a real vector $Y(L_i)\in H^1(T^2;\R)$ in the notation by Joyce {\rm\cite[Definition 6.2 and (72)]{J5}.}

We make the following hypothesis, which we verify shortly in our example:
\begin{hp}\l{H2}
There be $L_X\in\{L_1,L_2,L_3\}$ such that
\e\l{HS}\Im[H^1(L_X;\Z)\to H^1(T^2;\Z)]=\Im[H^1(X';\Z)\to H^1(T^2;\Z)].\e
This $L_X$ will be determined uniquely because $y(L_i),$ $i\in\{1,2,3\},$ are pairwise linearly-independent.
\end{hp}

Define $(N_1;\om_1,\Om_1;X_1)$ as in \S\ref{con} but with respect to the example of Theorem \ref{EB}; then:
\begin{thm}\l{P2}
Hypothesis $\ref{H2}$ holds with $X=X_1.$
\end{thm}
\begin{proof}
Let
$L_{X_1}:=\{(x_2,x_4,y)\in\C^3:|x_2|^2=|x_4|^2-1=\de|y|^2,\re x_2x_4y\in\R_{\ge0}\}$
where $(x_2,x_4,y)$ is as in the proof of Theorem \ref{NCBB}.
Then, $H_1(L_{X_1};\R)$ is generated by a 1-cycle on which $x_2$ is constant and $x_4$ is parametrized by $S^1.$
This 1-cycle corresponds to an $S^1$-orbit in Lemma \ref{hm} (b), which implies \eq{HS}.
\end{proof}

We denote by $\cM(\om)$ the moduli space of closed integral special Lagrangians in $(M;\om,\Om)$ as currents in the sense of geometric measure theory.
We regard $X$ as an element of $\cM(\om)$ as a multiplicity-one current in $M.$
The present author \cite{I} proves:
\begin{thm}\l{M}
The germ of $X$ in $\cM(\om)$ is homeomorphic to $(\R_{\ge0},0)$ where $\R_+$ corresponds to the desingularizations of $X$ by Joyce {\rm\cite[Theorem 10.4]{J5}} modelled upon $L_X,$ all diffeomorphic to $S^1\t S^2$ by a genus-one Heegaard splitting associated to $\eq{HS}.$
\end{thm}
\begin{cor*}
Theorem $\ref{M}$ applies to $(N_1';\om_1,\Om_1;X_1).$
\end{cor*}

We present two versions of Theorem \ref{M}, which are Theorems \ref{T1} and \ref{T2};
we apply Theorem \ref{T1} to the smoothings of $N_1,$ and Theorem \ref{T2} to the projective small resolutions of $N_1.$

Though the two versions require modification of the proof of Theorem \ref{M}, yet the methods are the same; which we therefore recall now in the case of Theorem \ref{M}.
There is a notion of energy such that every special Lagrangian in $(M;\om,\Om)$ close enough to $X$ is either of small energy or of large energy; and in the former case, isotopic to $X$ whereas, in the latter case, topologically a gluing of $X$ and $L_i$ for some $i\in\{1,2,3\}.$
On the other hand, since $H^1_c(X';\R)=0$ it follows that $X$ is rigid with respect to those deformations of $X$ by Joyce \cite{J2} which preserve the singularity of $X.$ Consequently, every nontrivial nearby special Lagrangian belongs to the latter. 
The gluing of $X$ and $L_i$ being Lagrangian, together with \eq{HS}, implies $L_i=L_X;$ whence one gets Theorem \ref{M}.

We come now to one of the two versions of Theorem \ref{M}.
We denote by $\io:X\to M$ the inclusion.
\begin{thm}\l{T1}
Suppose given $l\in\Z_{\ge0}$ and, on $M,$ an $l$-dimensional $C^\iy$-family of K\"ahler forms $\om^t,$ $t\in\R^l$ with $|t|$ small enough, $\om^0=\om,$ and $\io^*[\om^t]\=0;$
then, the germ of X in $\bigcup_t\cM(\om^t)$ with respect to the topology of currents in $M$ is homeomorphic to that of $0$ in $\R_{\ge0}\t\R^l$ where $\{0\}\t\R^l$ corresponds to the perturbations of $X$ by Joyce {\rm\cite[Corollary 5.8]{J5},} and $\R_+\t\R^l$ to their desingularizations by $L_X$ $($and accordingly, as in Theorem $\ref{M},$ all diffeomorphic to $S^1\t S^2).$ 
\end{thm}
\begin{proof}
The nearby special Lagrangians of small energy are the perturbations of $X;$ and those of large energy are their desingularizations.
\end{proof}

The other version of Theorem \ref{M}, which we explain now, is more complicated.
The properties of $y(L_i)$ imply that for every nonzero vector $v\in H^2(X;\R)=\R$ there exists a unique $L_v\in\{L_1,L_2,L_3\}$ such that the image in $\eq{ES}$ of $y(L_v)\in H^1(T^2;\R)$ is proportional to $v$ by a positive constant; in particular, $L_v$ and $L_X$ are mutually distinct with
\e\l{VW}H^1(T^2;\Z)=\Z y(L_v)+\Z y(L_X).\e
\begin{thm}\l{T2}
Suppose given on $M$ a one-dimensional $C^\iy$-family of K\"ahler forms $\om^t,$ $t\in\R$ and $|t|$ small enough, $\om^0=\om,$ and
\e\l{n0}0\ne v:=\frac d{dt}\io^*[\om^t]\Bigm|_{t=0}\in H^2(X;\R);\e
then, the germ of $X$ in $\bigcup_t\cM(\om^t)$ is homeomorphic to that of $0$ in $\R_{\ge0}$ where $\R_+$ corresponds to the desingularizations of $X$ by Joyce {\rm\cite[Theorem 10.5]{J5}} modelled upon $L_v,$ all diffeomorphic to $S^3$ by a genus-one Heegaard splitting associated to \eq{VW}.
\end{thm}
\begin{proof}
The hypothesis \eq{n0} implies $\io^*[\om^t]\ne0,$ $|t|$ small enough; in particular, no Lagrangian in $(M,\om^t)$ is isotopic to $X;$ and accordingly, every nontrivial nearby special Lagrangian is of large energy, topologically a gluing of $X$ and $L_i$ for some $i.$
We prove $L_i=L_v.$

There is, as Joyce \cite[Theorem 4.8]{J1} constructs, a family of Weinstein-like neighbourhoods of $X'$ on which, with $t$ fixed, $\om^t$ is the sum of
\iz
\item the pull-back under the projection $T^*X'\to X'$ of a compactly-supported 2-form $\be^t$ with $[\be^t]=\io^*[\om^t]\in H_c^2(X';\R)=H^2(X;\R);$ and
\item the canonical symplectic form on $T^*X'$ whose sign is such that it agrees, on the graph of a 1-form $\al$ on $X',$ with $-d\al.$
\iz
Recall that $X^t$ is, outside a neighbourhood of $x$ in $M,$ expressible as the graph of a 1-form $\al^t$ on a compact manifold-with-boundary $K\sb X'.$
Then, $X^t$ being Lagrangian in $(M,\om^t)$ implies $-d\al^t+\be^t=0;$ and in particular,
\e\l{al}[d\al^t]=[\be^t]=\io^*[\om^t]\in H^2(X;\R).\e
On the other hand, Joyce \cite[Proposition 7.6]{J1} proves that $L_i$ is, near infinity, expressible as the graph of a closed 1-form whose cohomology class is proportional to $y(L_i)$ by a $t$-independent positive constant.
Hence, recalling that $L_i$ is the blow-up of $X^t,$ one finds $\de^t>0$ and $c_1>0,$ the latter independent of $t,$ such that
\e\l{de}\de^t\{c_1y(L_i)+o(1)\}=[\al^t|_{\d K}] \e
where $o(1)$ converges to $0$ as $t\to0.$
On the other hand, since $\al^0=0$ it follows that $[\al^t|_{\d K}]=O(t).$
This combined with \eq{de} implies $\de^t=c_2t+o(1),$ $c_2$ a $t$-independent positive constant.
Then, by \eq{n0}, the image of $y(L_i)$ in \eq{ES} is $(c_1c_2)^{-1}v;$ which implies $L_i=L_v.$

Thus, every nontrivial nearby special Lagrangian is a gluing of $X$ and $L_v.$
\end{proof}

\subsection{Projective Small Resolutions of $N_1$}\l{23}
We return now to the example.
We denote by $D_0^{ij},$ $(i,j)\in\{1,4\}\t\{2,3\},$ the prime divisor on $N_0$ defined by $z_i=z_j=0,$ which implies in particular $D_0=D_0^{34};$ and denote by $D_1^{ij}$ the strict transform of $D_0^{ij}$ under $\pi:N_1\to N_0.$
\begin{lem*}
There exist projective small resolutions of $N_1.$
\end{lem*}
\begin{proof}
The successive blow-ups in any order of the three divisors $D_1^{ij},$ $(i,j)\ne(3,4),$ define projective small resolutions of $N_1.$
\end{proof}
\begin{lem}\l{B5}
Let $N_2\to N_1$ be a projective small resolution. Then, $\pi_1N_2=1;$ $b_2(N_2)=5;$ $\Om^3_{N_2}\cong\O_{N_2};$ and $H^1(N_2,\O_{N_2})=H^2(N_2,\O_{N_2})=0.$
\end{lem}
\begin{proof}
The embedding $N_0\sb\CP^4$ induces, by Lefschetz's hyperplane-section theorem, $\pi_1(N_0)=\pi_1(\CP^4)=1.$ Hence, recalling that $N_2\to N_1$ and $N_1\to N_0$ are both small morphisms, one finds, for every $k\in\{1,2,3\},$ $\pi_1N_k=1$ and $\Om^3_{N_k}\cong\O_{N_k}.$

In particular, the nonsingular 3-fold $N_2$ admits, according to Calabi--Yau and Bogomolov, a Riemannian metric of holonomy $SU_3.$ This implies in turn that $H^p(N_2,\O),$ $p\in\{1,2\},$ is isomorphic to the $SU_3$-invariant subspace of $\bigwedge^p\C^3;$ the latter, in fact, vanishes.

We prove $b_2(N_2)=5.$ Recall that Betti numbers are independent of projective small resolutions, a few of the flop-invariants (given by Koll\'ar \cite[Theorem 3.2.2]{Koll} for instance).
One may then suppose that $N_2$ is one of the successive blow-ups of the three divisors $D_1^{ij},$ $(i,j)\ne(3,4),$ and in particular that there is an embedding $N_2\sb N_1\t(\CP^1)^3.$
There are two other embeddings $N_0\sb\CP^4$ and $N_1\sb N_0\t\CP^1.$
These induce, by Lefschetz's hyperplane-section theorem, $b_2(N_0)=b_2(\CP^4)=1,$ $b_2(N_1)=b_2(N_0\t\CP^1)=2,$ and $b_2(N_2)=b_2(N_1\t(\CP^1)^3)=5.$
\end{proof}

Recall also from birational geometry that nodes are rational singularities; then, one finds from Lemma \ref{B5}:
\begin{cor}\l{NK}
$H^1(N_1,\O_{N_1})=H^2(N_1,\O_{N_1})=0.$
\end{cor}

Let $N_2\to N_1$ be again a projective small resolution.
We denote by $D_2^{ij},$ $(i,j)\in\{1,4\}\t\{2,3\},$ the strict transform of $D_1^{ij}$ under $N_2\to N_1.$
% and put $D_k:=D_k^{34},$ $k\in\{1,2\},$ extending $D_0=D_0^{34}.$

We denote by $H_0$ the hyperplane in $N_0\sb\CP^4$ defined by $z_0=0.$
We denote by $H_1$ the strict transform of $H_0$ under $\pi:N_1\to N_0,$ and by $H_2$ the strict transform of $H_1$ under $N_2\to N_1.$

We denote by $c_1(D)$ the first Chern class over $\R$ of a Cartier divisor $D.$
\begin{lem*}
%We introduce a notation for first Chern classes.
%Let $N$ be a projective nodal complex 3-fold, and $D$ a prime divisor on $N$ which does not intersect the nodes of $N;$ then, $D$ is a Cartier divisor on $N$ and accordingly induces an invertible sheaf over $N;$ whose first Chern class with real coefficients we denote by $c_1(D)\in H^2(N;\R).$
%It is, if $N$ nonsingular, the Poincar\'e dual therein to the 4-cycle $D.$
\begin{align}
\l{N0}H^2(N_0;\R)&=\R c_1(H_0)=\R[\om_0];\\
\l{N1}H^2(N_1;\R)&=\R c_1(H_1)+\R c_1(D_1^{34})=\R[\pi^*\om_0]+\R[\om_1];\text{ and }\\
\l{N2}H^2(N_2;\R)&=\R c_1(H_2)+\R\langle c_1(D_2^{ij})\rangle_{(i,j)\in\{1,4\}\t\{2,3\}}.
\end{align}
\end{lem*}
\begin{proof}
The dimensions of these three vector-spaces are given already by Lemma \ref{B5} and its proof.
On the other hand, the first Chern classes in these three lines are linearly independent, respectively.
Hence one gets the first equalities of the three lines.
The second equalities of the first two lines are obvious from the definitions of $\om_0$ and $\om_1.$
\end{proof}
%Let $\Nh$ be the blow-up of $D$ in $N,$ and $\Dh$ the strict transform of $D$ under it. Then, the embedding $\Nh\sb N\t\CP^1$ induces, by Lefschetz's hyperplane-section theorem, an isomorphism $H^2(\Nh;\R)=H^2(N\t\CP^1;\R);$ and the homomorphism $H^2(\CP^1;\R)\to H^2(\Nh;\R)$ induced by the projection $N\to\CP^1$ maps the Fubini--Study class to $-c_1(\Dh).$

We denote by $X_2$ the pre-image of $X_1$ under $N_2\to N_1;$ and by $\io_k:X_k\to N_k,$ $k\in\{0,1,2\},$ the inclusion.
% which is, outside the exceptional fibres of $\pi_2:N_2\to N_1,$ special Lagrangian in the ordinary sense with respect to $(\pi_2^*\om_1,\pi_2^*\Om_1).$
\begin{lem}\l{SC}{\bf(a)}
The homomorphism $\io_k^*:H^2(N_k;\R)\to H^2(X_k;\R),$ $k\in\{1,2\},$ vanishes; and {\bf(b)}
the homomorphism $\io_2^*:H^2(N_2;\R)\to H^2(X_2;\R)$ assigns to every K\"ahler class on $N_2$ a nonzero element.
\end{lem}
%Let $\Si_0$ be a 2-cycle in $X_0$ which generates $H_2(X_0;\R),$ and $\Si_1$ a 2-cycle of $X_1$ which generates $H_2(X_1;\R).$
%Recalling again that $X_1$ does not intersect the nodes of $N_1,$ we may well define $\Si_2:=\pi_2^*\Si_1,$ which generates $H_2(X_2;\R).$
%
%$\io_{k*}[\Si_k]=0,$ $k\in\{0,1\};$ and $\io_{2*}[\Si_2]\cap c_1(H_2)=\io_{2*}[\Si_2]\cap c_1(D_2)=\io_{2*}[\Si_2]\cap c_1(D_2^{24})=0$ whereas $\io_{2*}[\Si_2]\cap c_1(D_2^{12}),$ $\io_{2*}[\Si_2]\cap c_1(D_2^{13}),$ and $\io_{2*}[\Si_2]\cap[\om_2]$ are all nonzero, $\om_2$ an arbitrary K\"ahler form on $N_2.$
\begin{proof}[Proof of {\rm(a)}]
It follows from \eq{N0}, \eq{N1} and the fact that $X_k,$ $k\in\{1,2\},$ is Lagrangian in $(N_k,\om_k).$
\end{proof}
\begin{proof}[Proof of {\rm(b)}]
The image of $\io_{2*}:H_2(X_2;\R)\to H_2(N_2;\R)$ is generated by a 2-sphere $\Si$ in a Lagrangian smoothing of $X_2$ where $\Si$ is, near the singular point of $X_2,$ defined by $x_2=\bar y$ and $x_4=(|y|^2+1)^{1/2},$ with $(x_2,x_4,y)$ as in the proof of Theorem \ref{P2}. In particular, $x_4\ne0$ on $\Si;$ accordingly, $\Si$ intersects neither $D_2^{34}$ nor $D_2^{24};$ and consequently, $\Si\cap c_1(D_2^{ij})=0,$ $(i,j)=(3,4)\text{ or }(2,4).$
On the other hand, $D_2^{12}\cap\Si=\{(0,0,0,1;0,1)\in Q_1\};$ at which point the tangent spaces to $D_2^{12}$ and $\Si$ are $\{(0,x_4,y):x_4,y\in\C\}$ and $\{(\bar y,0,y):y\in\C\},$ respectively, mutually transverse.
Likewise, $D_2^{13}\cap\Si$ is the same one-point set; at which the tangent spaces to $D_2^{13}$ and $\Si$ are $\{(x_2,x_4,0):x_2,x_4\in\C\}$ and $\{(\bar y,0,y):y\in\C\},$ respectively, also mutually transverse but oppositely in the sense that $[\Si]\cap c_1(D_2^{12})=-[\Si]\cap c_1(D_2^{13})=\pm1$ according to the orientation of $\Si.$

On $N_2$ every K\"ahler class $[\om_2]$ is as in \eq{N2} an $\R$-linear combination of the first Chern classes of the five divisors; whose coefficient of a divisor $D$ we denote by $t(D)\in\R.$
Integrating $\om_2$ over the exceptional fibres over the nodes at which exactly two of $\{z_1,z_2,z_3,z_4\}$ vanish, one finds $t(D_2^{ij})>0,$ $\{i,j\}\in\{1,4\}\t\{2,3\}.$
%One has also $t(H_0)>0$ which is however unnecessary at the moment.
Let $E$ be the exceptional fibre of a node of $N_1$ at which $z_1=z_2=z_3=0$ and $z_4\ne0.$
Then, near $E,$ the blow-ups of $D_2^{12}$ and $D_2^{13}$ are related by a flop; in particular, $[E]\cap c_1(D_2^{12})=-[E]\cap c_1(D_2^{13}).$ Hence one finds, integrating $\om_2$ over $E,$ $t(D_2^{12})\ne t(D_2^{13});$
and consequently, $[\Si]\cap[\om_2]=\pm[t(D_2^{12})-t(D_2^{13})]\ne0.$
\end{proof}

We denote by $\Om_2$ the pull-back of $\Om_1$ under $N_2\to N_1,$ which defines a holomorphic volume-form on $N_2;$ and by $M_2$ the complement in $N_2$ to the exceptional fibres over $N_1.$
Since $X_1$ does not intersect the nodes of $N_1$ it follows that $\io_2^*\om_1$ is positive definite on $M_2$ and that $X_2$ is special Lagrangian in $(M_2;\io_2^*\om_1,\Om_2).$
Lemma \ref{SC} (b) implies:
\begin{cor*}
Suppose given an arbitrary one-dimensional smooth family $\om_2^t,$ $t\in\R,$ with $\om_2^0=\io_2^*\om_1,$ of closed real $(1,1)$-forms on $N_2;$ suppose $\om_2^t,$ $t>0,$ positive definite $($and accordingly K\"ahlerian$);$ and suppose nonzero the vector $w:=d[\om_2^t]/dt|_{t=0}\in H^2(N_2;\R).$
Then, $\io_2^*w\ne0;$ and consequently, Theorem {\rm\ref{T2}} applies to $(M_2;\om_2^t,\Om_2;X_2).$
\end{cor*}

\subsection{Smoothings of $N_1$}\l{24}
We recall the following general fact:
\begin{thm}\l{AG}
Let $N$ be a projective nodal complex $3$-fold with $\Om_N^3\cong\O_N$ and $H^1(N,\O_N)=H^2(N,\O_N)=0;$ then:
\item{\bf(a)}
$N$ admits an unobstructed universal deformation $\varpi:\mathfrak{N}\to\De.$ It is projective over $\De;$ and in fact, every invertible sheaf over $N$ extends to $\mathfrak{N}$ with $\De$ made smaller if need be. Moreover, $R^0\varpi_*\Om^3_{\mathfrak{N}/\De}$ is of constant rank one over $\De.$
\item{\bf(b)}
$N$ be smoothable if and only if there exist a small resolution $\pi:M\to N$ and, for all the nodes $n\in N,$ some nonzero real numbers $c_n$ with $\sum c_n[\pi^{-1}(n)]=0\in H_2(M;\R);$ the latter condition is in fact independent of the choice of the small resolution $\pi:M\to N.$
\end{thm}
\begin{proof}
Kawamata \cite[Corollary 8.6]{Kaw1} proves $H^0(N,T_N)\cong H^1(N,\O_N),$ which vanishes by hypothesis; consequently, $N$ admits a universal deformation. It is unobstructed according to Kawamata \cite[Theorem 4]{Kaw2} and Tian \cite[Theorem 0.2]{T}. The next assertion in (a) follows since $H^2(N,\O_N),$ which vanishes by hypothesis, is the obstruction to extending an invertible sheaf over $N.$
The last assertion of (a) follows since $H^1(N,\O_N)=0.$
% which vanishes by hypothesis, is the obstruction to extending a section of $\O_N\cong\Om_N^3.$

Part (b) follows from the result by Friedman \cite[Proposition 8.7]{Fr} together with the unobstructedness in (a).
\end{proof}
Corollary \ref{NK} implies that $N_1$ satisfies the hypothesis of Theorem \ref{AG}; and in fact:
\begin{cor*}
$N_1$ is smoothable.
\end{cor*}
\begin{proof}
In $N_0$ or $N_1$ the nodes of type $I\sb\{1,2,3,4\}$ shall mean those with $z_i=0$ if and only if $i\in I.$
The divisor $D_0=D_0^{34}$ contains sixteen nodes in $N_0$: the one of type $\{1,2,3,4\}$, the three of type $\{1,3,4\}$, the three of type $\{2,3,4\},$ and the nine of type $\{3,4\}.$
The number of the nodes of $N_1$ is $49-16=33$ which are counted also as follows: the three of type $\{1,2,3\},$ the three of the type $\{1,2,4\},$ the nine of type $\{1,2\},$ the nine of type $\{1,3\},$ and the nine of type $\{2,4\}.$
By \eq{N2}, the homology class of an exceptional fibre of the small resolution $N_2\to N_1$ is expressible as a vector of the three numbers of intersections with $D_2^{12},$ $D_2^{13}$ and $D_2^{24}$ in this order.
This vector depends only upon the type of the node over which the exceptional fibre is; and those vectors corresponding to type $\{1,2,3\},$ $\{1,2,4\},$ $\{1,2\},$ $\{1,3\},$ and $\{2,4\}$ are $v_1:=(1,1,0),$ $v_2:=(1,0,1),$ $v_3:=(1,0,0),$ $v_4:=(0,1,0),$ and $v_5:=(0,0,1)$ respectively.
There is a relation $-v_1-v_2+2v_3+v_4+v_5=0,$ all the coefficients nonzero; which implies the criterion of Theorem \ref{AG} (b).
\end{proof}

%We treat in what follows, for the brevity of notations, an arbitrary but one-dimensional family of smoothings of $N_1.$ We denote it by $\varpi:\mathfrak{N}\to\De,$ where $\De$ is a neighbourhood of 0 in $\C;$ identify $\varpi^{-1}(0)$ with $ N_1;$ and shrink $\mathfrak{N},\De$ wherever need be.

The affine nodal $3$-fold $Q_0,$ defined in $\C^4$ by $x_1x_4-x_2x_3=0,$ admits a semi-universal deformation $\mathfrak{U}\to \De$ where $\De$ is a neighbourhood of 0 in $\C,$ $\mathfrak{U}$ a nonsingular hypersurface in $\C^4\t\De$ defined by $x_1x_4-x_2x_3=s,$ and the map $\mathfrak{U}\to \De$ defined by the co\"ordinate $s.$
We shrink $\De$ wherever need be in what follows.

Let $\varpi:\mathfrak{N}\to\De$ be now a one-dimensional family of smoothings of $N_1,$ which we identify with $\varpi^{-1}(0),$ such that for every node $n\in N_1$ there exists a biholomorphism of germs $(\mathfrak{N},n)\cong(\mathfrak{U},0)$ which commutes with the two projections onto $\De;$ in particular, $\mathfrak{N}$ is nonsingular.

By a {\it $C^\iy$-family of K\"ahler forms on $M^s:=\varpi^{-1}(s)$} we mean a real 2-form on $\mathfrak{N}$ which is:
\iz
\item
outside the nodes of $ N_1,$ locally of the form $\sum_{a,b=1}^3ig_{a\bar b}(x_1,x_2,x_3;s)dx_a\w d\bar{x}_b$ with $g_{a\bar b}=\ov{g_{b\bar a}}$ a $C^\iy$-function, and $(x_1,x_2,x_3;s)$ a local co\"ordinate-system in $\mathfrak{N}$ with $s=\varpi;$ and
\item
near every node $n\in N_1,$ of the form $\sum_{a,b=1}^4ig_{a\bar b}(x_1,x_2,x_3,x_4;s)dx_a\w d\bar{x}_b$ with $g_{a\bar b}=\ov{g_{b\bar a}}$ a $C^\iy$-function, and $(x_1,x_2,x_3,x_4;s)$ the co\"ordinates of $\mathfrak{U}$ under an identification $(\mathfrak{N},n)\cong(\mathfrak{U},0).$
\iz
\begin{lem*}
$\om_1$ extends to a $C^\iy$-family of K\"ahler forms on $M^s.$
\end{lem*}
\begin{proof}
We introduce a notion of Hermitian metrics on invertible sheaves on $N_1.$
Let $\{U_i\}$ be a covering of $ N_1,$ and $\{f_{ij}\}\in C^1( N_1,\O^*)$ a \v Cech cocycle with respect to $\{U_i\}$ which defines of course an invertible sheaf over $ N_1;$ by a {\it Hermitian metric} on which, then, we mean a system of $C^\iy$-functions $h_i:U_i\to\R_+$ with $h_i=|f_{ij}|^2h_j$ where, if $U_i$ contains a node, a $C^\iy$-function on $U_i$ means the restriction to $U_i$ of some $C^\iy$-function near 0 in $\C^4$ under a holomorphism-of-germs $(U_i,n)\cong(Q_0,0)\sb(\C^4,0).$
We denote by $c_1\{h_i\}$ the {\it curvature-form} of the Hermitian metric $\{h_i\};$ which is, defined in the same manner as in the nonsingular case, locally of the form $\sqrt{-1}\d\db\log h_i.$

The K\"ahler form $\om_1$ is then of the form $c_1(h')-\ep^2c_1(h'')$ where $h',h''$ are Hermitian metrics on the invertible sheaves associated to the Cartier divisors $H_1,D_1$ respectively.
These two invertible sheaves, by Theorem \ref{AG} (b), extend to $\mathfrak{N};$ to which, in turn, the two Hermitian metrics also extend with a partition of unity.
Taking their fibrewise curvature-forms and their linear combination with the factor $\ep^2$ one gets a required extension of $\om_1.$
\end{proof}

\begin{cor*}\l{persist}
{\bf(a)}
For every smooth family $\om^s,$ with $\om^0=\om_1,$ of K\"ahler forms on $M^s$ there exist a $C^\iy$-section $s\mapsto\Om^s$ of $R^0\varpi_*\Om^3_{\mathfrak{N}/\De}\cong\O_\De,$ with $\Om^0=\Om_1,$ and a continuous family $X^s,$ with $X^0=X_1,$ of compact special Lagrangians in $(M^s;\om^s,\Om^s)$ with one-point singularities modelled upon Harvey--Lawson's $T^2$-cone.
\item{\bf(b)}
For any such $(\om^s,\Om^s;X_1^s)$ as in {\rm(a)} and for each $s\ne0,$ the homomorphism $\io_1^*:H^2(M^s;\R)\to H^2(X_1^s;\R)$ is the zero map; and consequently, Theorem $\ref{T1}$ applies to $(M^s;\om^s,\Om^s;X_1^s)$ together with any $C^\iy$-family $\om^{s,t},$ $t\in\R^l$ and $|t|$ small enough, with $\om^{s,0}=\om^s,$ of K\"ahler forms on $M^s.$
\end{cor*}

\begin{proof}
Since $X_1\sb N_1'$ it follows that the inclusion $X_1\sb N_1$ extends to a $C^\iy$-family of embeddings $X_1\sb M^s$ which we use.
One finds then, extending $\Om_1$ to a section of $R^0\varpi_*\Om^3_{\mathfrak{N}/\De}\cong\O_\De$ and applying to it a suitable scalar multiplication, a $C^\iy$-section $s\mapsto\Om^s$ of $R^0\varpi_*\Om^3_{\mathfrak{N}/\De}$ with $[X_1]\cap[\Im\Om^s]\=0.$
With this $\Om^s$ the result by Joyce \cite[Corollary 5.8]{J5} yields $X^s$ as in (a).

Part {\rm(b)} follows from Lemma \ref{SC} (a) with $k=1$ and from the natural isomorphism $H^2(M^s;\R)\cong H^2(N_1;\R)$ induced by collapsing the 3-spheres associated to the nodes of $N_1.$
\end{proof}

ShanghaiTech University, 393 Middle Huxia Road, Pudong Shanghai, 201210 China

e-mail address: YosukeImagi@shanghaitech.edu.cn

\end{document}